\documentclass[%
pre,
longbibliography,
reprint,
superscriptaddress,
showpacs,
showkeys,
 amsmath,amssymb,amsthm
]{revtex4-1}

\usepackage{graphicx}% Include figure files
\usepackage{dcolumn}% Align table columns on decimal point
\usepackage{bm,bbm}% bold math
\usepackage{mathtools}
\usepackage[normalem]{ulem}

\renewcommand{\vec}[1]{\mathbf{#1}}
\renewcommand{\d}{{\rm d}}
\newcommand{\e}{{\rm e}}

\newcommand{\FD}[2]{\frac{\d #1}{\d #2}}

\DeclareSymbolFont{AMSb}{U}{msb}{m}{n}
\DeclareMathSymbol{\Nset}{\mathalpha}{AMSb}{"4E}
\DeclareMathSymbol{\Zset}{\mathalpha}{AMSb}{"5A}
\DeclareMathSymbol{\Rset}{\mathalpha}{AMSb}{"52}
\DeclareMathSymbol{\Qset}{\mathalpha}{AMSb}{"51}
\DeclareMathSymbol{\Cset}{\mathalpha}{AMSb}{"43}

\begin{document}

%\preprint{APS/123-QED}

\title{Clusters in nonsmooth oscillator networks}% Force line breaks with \\

\author{Rachel Nicks}
\email{rachel.nicks@nottingham.ac.uk}
\affiliation{Centre for Mathematical Medicine and Biology, School of Mathematical Sciences,\\
University of Nottingham,NG7 2RD, UK.
}

\author{Lucie Chambon}
\email{lucie.chambon@inria.fr}
\affiliation{
Centre de recherche INRIA Sophia-Antipolis M\'editerran\'ee, Borel building\\
2004, route des Lucioles -- BP 93\\
06 902 Sophia Antipolis Cedex., France.
}

\author{Stephen Coombes}
\email{stephen.coombes@nottingham.ac.uk}
\affiliation{Centre for Mathematical Medicine and Biology, School of Mathematical Sciences,\\
University of Nottingham,NG7 2RD, UK.
}

\date{\today}

\begin{abstract}
For coupled oscillator networks with Laplacian coupling the master stability function (MSF) has proven a particularly powerful tool for assessing the stability of the synchronous state.  Using tools from group theory this approach has recently been extended to treat more general cluster states.  However, the MSF and its generalisations require the determination of a set of Floquet multipliers from variational equations obtained by linearisation around a periodic orbit.  Since closed form solutions for periodic orbits are invariably hard to come by the framework is often explored using numerical techniques.  Here we show that further insight into network dynamics can be obtained by focusing on piecewise linear (PWL) oscillator models.  Not only do these allow for the explicit construction of periodic orbits, their variational analysis can also be explicitly performed.  The price for adopting such nonsmooth systems is that many of the notions from smooth dynamical systems, and in particular linear stability, need to be modified to take into account possible jumps in the components of Jacobians.  This is naturally accommodated with the use of \textit{saltation} matrices.  By augmenting the variational approach for studying smooth dynamical systems with such matrices we show that, for a wide variety of networks that have been used as models of biological systems, cluster states can be explicitly investigated.  By way of illustration we analyse an integrate-and-fire network model with event-driven synaptic coupling as well as a diffusively coupled network built from planar PWL nodes, including a reduction of the popular Morris--Lecar neuron model.  We use these examples to emphasise that the stability of network cluster states can depend as much on the choice of single node dynamics as it does on the form of network structural connectivity.  Importantly the procedure that we present here, for understanding cluster synchronisation in networks, is valid for a wide variety of systems in biology, physics, and engineering that can be described by PWL oscillators.
\end{abstract}

%Synchronization; coupled oscillators
%Bifurcation nonlinear dynamics,
%Group theory mathematics,


\maketitle

\section{\label{sec:level1}Introduction}

The study of synchrony in coupled oscillator networks has a very long history, dating all the way back to the work of Huygens on two interacting pendulum clocks.  Since his observations about the emergence of ``an odd kind of sympathy" there have been countless other examples of synchrony discussed in the natural sciences and engineering ranging from the dynamics of populations of flashing fireflies to those of coupled Josephson junctions.  For an excellent review see Arenas \textit{et al.} \cite{Arenas2008} and the recent focus issue on patterns of network synchronisation \cite{Abrams2016}.  However, perfect global synchronisation is just one of many states expected to emerge in structured oscillator networks.  Indeed instabilities of the synchronous state are generically expected to give rise to cluster states, in which sub-populations may synchronise, though not with each other.  This class of solutions has been relatively well explored for phase-oscillator networks, as in the work by Brown \textit{et al}. \cite{Brown03} and Ashwin \textit{et al}. \cite{Ashwin07}, though less so for networks of limit-cycle oscillators.  For this more general scenario Golubitsky and Stewart have made use of the Hopf bifurcation with symmetry to understand cluster states in ring networks, and provided a number of case studies of systems with nearest-neighbour coupling \cite{Golubitsky1986}.  Pogromsky and collaborators have also exploited symmetry under permutation of a given network of dynamical systems coupled through diffusion to determine the stability of cluster states, with the aid of an appropriate Lyapunov function \cite{Pogromsky2002,Pogromsky2008}.  For the case of diffusively coupled oscillators, Belykh \textit{et al}. have further provided a complete classification of cluster states, together with conditions that can be used to determine the coexistence of stable cluster states with an unstable synchronous state \cite{Belykh2000}.  More recently Pecora \textit{et al}. \cite{Pecora2014}  and Sorrentino \textit{et al}. \cite{Sorrentino2016}, have extended the master stability function (MSF) approach of Pecora and Carroll \cite{Pecora1998} making extensive use of tools from computational group theory.  This work paves the way for a systematic study of cluster states and their bifurcations for a very broad class of networks with state-dependent interactions, including arrays of electrochemical oscillators \cite{Kiss2007}, gene networks \cite{Zhang2009}, and the brain \cite{Zemonova2006}.  Moreover, it will likely play a key role in addressing the growing interest in dynamics on networks \cite{Porter2016}.

At heart the MSF approach exploits the fact that the Floquet theory for the stability of a synchronous network state decouples into a set of lower dimensional Floquet problems.  By studying just one of these lower dimensional Floquet problems, as a function of a single complex parameter, the linear stability of the synchronous network state can be determined.
Indeed the approach is sufficiently general that, for \textit{Laplacian} coupling of identical nodes (so that the network synchronous state is guaranteed to exist as long as a single node oscillates), then the effect of different network choices is easily explored.  Namely for a given choice of node dynamics the MSF maps a complex number to the largest Floquet exponent of the reduced problem, and by choosing the complex number to run over all possible eigenvalues of the unnormalised graph Laplacian of the network connectivity then stability is guaranteed if the MSF is always negative.  The challenge of exploring different networks thus reduces to determining the eigenstructure of the graph Laplacian, with the caveat that the Floquet problem can be solved.  However, for a general nonlinear system (without special structure) it is unlikely that this problem can be solved analytically.  Thus it would be highly complementary to the MSF formalism, and its extensions to treat cluster states, if this barrier could be reduced or removed entirely.  It is this issue that we address in this paper.  We do so by restricting the choice of node dynamics to that of piecewise linear (PWL) systems.  Although at first sight this may appear overly restrictive, there has been an appreciation for some time in the mathematical community of the benefits of studying caricatures of complex systems built from PWL and possibly discontinuous dynamical systems \cite{diBernardo2008}.
Indeed, there is a now growing perspective in the applied dynamical systems community that piecewise models are highly suitable for modern applications in science and can complement the smooth dynamical systems approach that has dominated to date \cite{Glendinning2015}.  It has recently been shown how to extend the MSF formalism for synchrony to PWL models with state-dependent coupling \cite{Coombes2016}.  Here we go further and show how to treat discontinuous systems of integrate-and-fire (IF) type, and how to analyse the existence and stability of more general cluster states in a broad class of PWL network models with both state and time-dependent interactions.

In \S \ref{sec:MSF} we give a brief recapitulation of the MSF for smooth dynamical systems as well as its extension to networks of nodes built from PWL oscillators.
In illustration of the utility of this combination of MSF and nonsmooth systems we show how to analyse synchrony in PWL IF networks with \textit{balanced} synaptic coupling.  Next in \S \ref{sec:Cluster} we survey the techniques from computational group theory that allow the discovery of cluster states from the topology of the network.  For PWL systems we then  show, in  \S \ref{sec:Nonsmooth}, how such orbits can be explicitly determined.  Similarly the variational equations that determine the stability of a cluster state are equally tractable.  By exploiting the PWL nature of the network dynamics we are able to determine the stability of cluster states without recourse to numerically evolving the variational equations.  Our approach also facilitates an explicit bifurcation analysis, which we illustrate for a simple five-node network with a PWL analogue of a Hopf normal form.
In \S \ref{sec:Applications} we further explore a node dynamics that exhibits a PWL analogue of a homoclinic bifurcation, as well as a PWL IF model.   We use these examples to highlight that the stability of network cluster states can depend as much on the choice of single node dynamics as it does on the form of network structural connectivity.  Finally in \S \ref{sec:Discussion} we discuss natural extensions of the work in this paper.

\section{\label{sec:MSF}The Master Stability Function:  a recapitulation}

The MSF allows the determination of the linear stability of the synchronous state in a quite large class of smooth networks of identical nodes.
To describe the MSF formalism it is convenient to consider $N$ nodes (oscillators) and let $\vec{x}_i \in \Rset^m$ be the $m$-dimensional vector of dynamical variables of the $i$th node with isolated (uncoupled) dynamics $\dot{\vec{x}}_i = \vec{F}(\vec{x}_i)$, with $i=1,\ldots, N$.  The output for each oscillator is described by a vector function $\vec{H} \in \Rset^m$.
For a given coupling matrix with components $W_{ij}$ and a global coupling strength $\sigma$ the network dynamics, to which the MSF formalism applies, is specified by
\begin{align}
\FD{}{t}{\vec{x}}_i &= \vec{F}(\vec{x}_i) + \sigma \sum_{j=1}^N W_{ij} \left [ \vec{H}(\vec{x}_j)-\vec{H}(\vec{x}_i) \right ] ,\nonumber \\
&\equiv \vec{F}(\vec{x}_i) - \sigma \sum_{j=1}^N \mathcal{G}_{ij} \vec{H}(\vec{x}_j) .
\label{network1}
\end{align}
Here the matrix $\mathcal{G}$ with components $\mathcal{G}_{ij}$ has the unnormalised graph Laplacian structure $\mathcal{G}_{ij} = -W_{ij} + \delta_{ij} \sum_{k}W_{ik}$.
The $N-1$ constraints $\vec{x}_1(t)=\vec{x}_2(t)=\ldots =\vec{x}_N(t) = \vec{s}(t)$ define the (invariant) synchronisation manifold, with $\vec{s}(t)$ a solution of the uncoupled system, namely $\dot{\vec{s}}=\vec{F}(\vec{s})$.  To assess the stability of this state a linear stability analysis is performed by expanding a solution as $\vec{x}_i(t) = \vec{s}(t) +\delta \vec{x}_i(t)$ to obtain the variational equation
\begin{equation}
\FD{}{t}\delta {\vec{x}_i} = D \vec{F}(\vec{s}) \delta \vec{x}_i - \sigma D \vec{H}(\vec{s}) \sum_{j=1}^N \mathcal{G}_{ij} \delta \vec{x}_j .
\nonumber
\end{equation}
Here $D \vec{F}(\vec{s})$ and $D \vec{H}(\vec{s})$ denote the Jacobian of $\vec{F}(\vec{s})$ and $\vec{H}(\vec{s})$ around the synchronous solution respectively.
The variational equation has a block form that can be simplified by projecting $\delta \vec{x}$ into the eigenspace spanned by the (right) eigenvectors of the matrix $\mathcal{G}$.  If we organise these column eigenvectors into a matrix $P$ then $\mathcal{G}P = P \Lambda$, with $\Lambda=\operatorname{diag}(\lambda_1,\ldots,\lambda_N)$, where $\lambda_l$ are the corresponding eigenvalues for $l=1,\dots,N$.  If we collect the perturbations in a vector $\vec{U}=(\delta \vec{x}_1, \ldots, \delta \vec{x}_N) \in \Rset^{Nm}$, and introduce a new variable $\vec{V}$ according to the linear transformation $\vec{V} = (P \otimes I_m)^{-1} \vec{U}$, then we have that
\begin{equation}
\FD{}{t}{\vec{V}} = (I_N \otimes D\vec{F}(\vec{s})) \vec{V} - \sigma(\Lambda \otimes D\vec{H}(\vec{s})) \vec{V} .
\label{Vdot}
\end{equation}
Here the symbol $\otimes$ denotes the tensor (or Kronecker) product for matrices, and $I_N$ is the $N \times N$ identity matrix.  Thus we have a set of $N$ decoupled eigenmodes in the block form
\begin{equation}
\FD{}{t}{\vec{\xi}_l} = \left [ D \vec{F}(\vec{s}) - \sigma \lambda_l D \vec{H}(\vec{s}) \right ] \vec{\xi}_l, \qquad l=1,\ldots,N,
\nonumber
\end{equation}
where $\vec{\xi}_l$ is the $l$th (right) eigenmode associated with the eigenvalue $\lambda_l$ of $\mathcal{G}$ (and $D \vec{F}(\vec{s})$ and $D \vec{H}(\vec{s})$ are independent of the block label $l$).  Since $\sum_j \mathcal{G}_{ij}=0$ there is always a zero eigenvalue, say $\lambda_1=0$, with corresponding eigenvector $(1,1,\ldots,1)$, describing a perturbation parallel to the synchronisation manifold.  The other $N-1$ transverse eigenmodes must damp out for synchrony to be stable.  For a general matrix $\mathcal{G}$ the eigenvalues $\lambda_l$ may be complex, which brings us to consideration of the system
\begin{equation}
\FD{}{t}{\vec{\xi}} = \left [ D \vec{F}(\vec{s}) - \beta D \vec{H}(\vec{s}) \right ] \vec{\xi}, \qquad \beta = \sigma \lambda_l \in \Cset  ,
\label{variational}
\end{equation}
where $\xi \in \Cset^m$.
The MSF is defined as the function which maps the complex number $\beta$  to the largest Floquet exponent of (\ref{variational}).
The synchronous state of the system of coupled oscillators is stable if the MSF is negative at  $\beta=\sigma \lambda_l$  where  $\lambda_l$  ranges over the eigenvalues of the matrix $\mathcal{G}$ (excluding $\lambda_1=0$).

We note that the \textit{Laplacian} form of coupling in (\ref{network1}) guarantees the existence of the synchronous state as long as a single uncoupled node has a periodic solution.  However, other forms of coupling, and in particular those described by adjacency matrices or their weighted counterparts are also common.  In this case (\ref{network1}) would be replaced by the system
\begin{equation}
\FD{}{t}{\vec{x}}_i = \vec{F}(\vec{x}_i) + \sigma \sum_{j=1}^N W_{ij} \vec{H}(\vec{x}_j) .
\label{network2}
\end{equation}
The synchronous solution (with an orbit inherited from the uncoupled node) is not a generic solution of (\ref{network2}) unless there is a further constraint placed on the coupling matrix, namely that $\sum_{j=1}^N W_{ij} = \text{constant}$ for all $i$.  If this constant value is $0$, then we shall say that the system is \textit{balanced}.  This scenario is ubiquitous in the modelling of many spiking neural networks \cite{VanVreeswijk1996}, and is relevant to understanding coding \cite{Deneve2016}, memory \cite{Roudi2007}, noise-robust neuronal selectivity \cite{Rubin2017}, and brain idling \cite{Doiron2014}.
Interestingly many spiking neural networks are modelled using integrate-and-fire (IF) neurons which are \textit{discontinuous} dynamical systems owing to the hard reset of the voltage state variable upon reaching a spiking threshold \cite{Coombes2011a}.  A naive application of the MSF formalism described above would lead to incorrect results since it assumes that the underlying dynamics is smooth.  Fortunately it is relatively straightforward to augment the MSF approach to handle firing discontinuities in IF networks using ideas from the study of impact oscillators.  This was originally done by Coombes to give meaning to Liapunov exponents for linear IF neurons \cite{Coombes99a}, and more recently by Ladenbauer \textit{et al}. \cite{Ladenbauer2013} for constructing the MSF in synaptically coupled networks of nonlinear (adaptive exponential) IF neurons. However, this latter work requires substantial numerical simulation as the periodic orbit for the synchronous solution is not available in closed form.  The perspective in this paper is to choose PWL caricatures of node dynamics to overcome this last barrier.  A case in point is the caricature of the adaptive exponential IF model developed in \cite{Naud2008},
which has a PWL subthreshold dynamics and an adaptive jump upon reset.  Since this is an exemplar PWL system we will describe it here in more detail to set the scene for the extension of the MSF formalism to cover nonsmooth systems in general and not just IF type models.  For a further perspective on the use of techniques from nonsmooth systems in neuroscience see \cite{Coombes2011a}.

\subsection{\label{sec:MSF-PWL}MSF for PWL IF networks with synaptic coupling}

There are now many variants on nonlinear integrate-and-fire neurons that are able to fit the spike trains of real neurons, such as those due to Gr\"obler \textit{et al}. \cite{Grobler98}, Izhikevich \cite{Izhikevich03}, and Badel \textit{et al}. \cite{Badel08}.  The planar adaptive exponential model has been particularly successful at fitting data from cortical fast spiking interneurons and regular spiking pyramidal neurons \cite{Badel08}.  Importantly it admits to a useful PWL reduction whereby its voltage nullcline is replaced by a PWL function \cite{Naud2008}.  A similar approach is taken in \cite{Coombes09a} to develop an analytically tractable PWL IF model that can support both tonic (repetitive firing with a constant interspike interval) and bursting states (with regular alternations between short and long interspike intervals).  Here we subsume both models within a general description of the form $\dot{\vec{z}}=\vec{F}(\vec{z})$, $\vec{z} \in \Rset^m$ with
\begin{equation}
\vec{F}(\vec{z}) = \begin{cases}
A_L \vec{z} + \vec{c}_L, & \vec{z} \in S_L \\
A_R \vec{z} + \vec{c}_R, & \vec{z} \in S_R
\end{cases},
\label{zdot}
\end{equation}
with $A_{L,R} \in \Rset^{m \times m}$, $\vec{c}_{L,R} \in \Rset^m$, and $\Rset^m = S_L \cup \Sigma \cup S_R$.
For simplicity we restrict attention to the situation that the phase-space of the model can be broken into two regions $S_{L,R}$ separated by a boundary given implicitly by the zero of an indicator function $h$, such that $\Sigma=\{ \vec{z} ~|~h(\vec{z})=0\}$.  It is a simple matter to break the phase-space up into more regions, to incorporate further \textit{switching} manifolds, so we shall stick to describing the simplest situation (though will give examples in \S \ref{sec:Applications} of systems with two switching manifolds).  Should the dynamics on the switching manifold not be continuous with that in regions $S_L$ and $S_R$ then it is common to adopt the Filippov convex method for dynamics on $\Sigma$ \cite{Filippov1998}.

\begin{figure}
\begin{center}
\includegraphics[width=8cm]{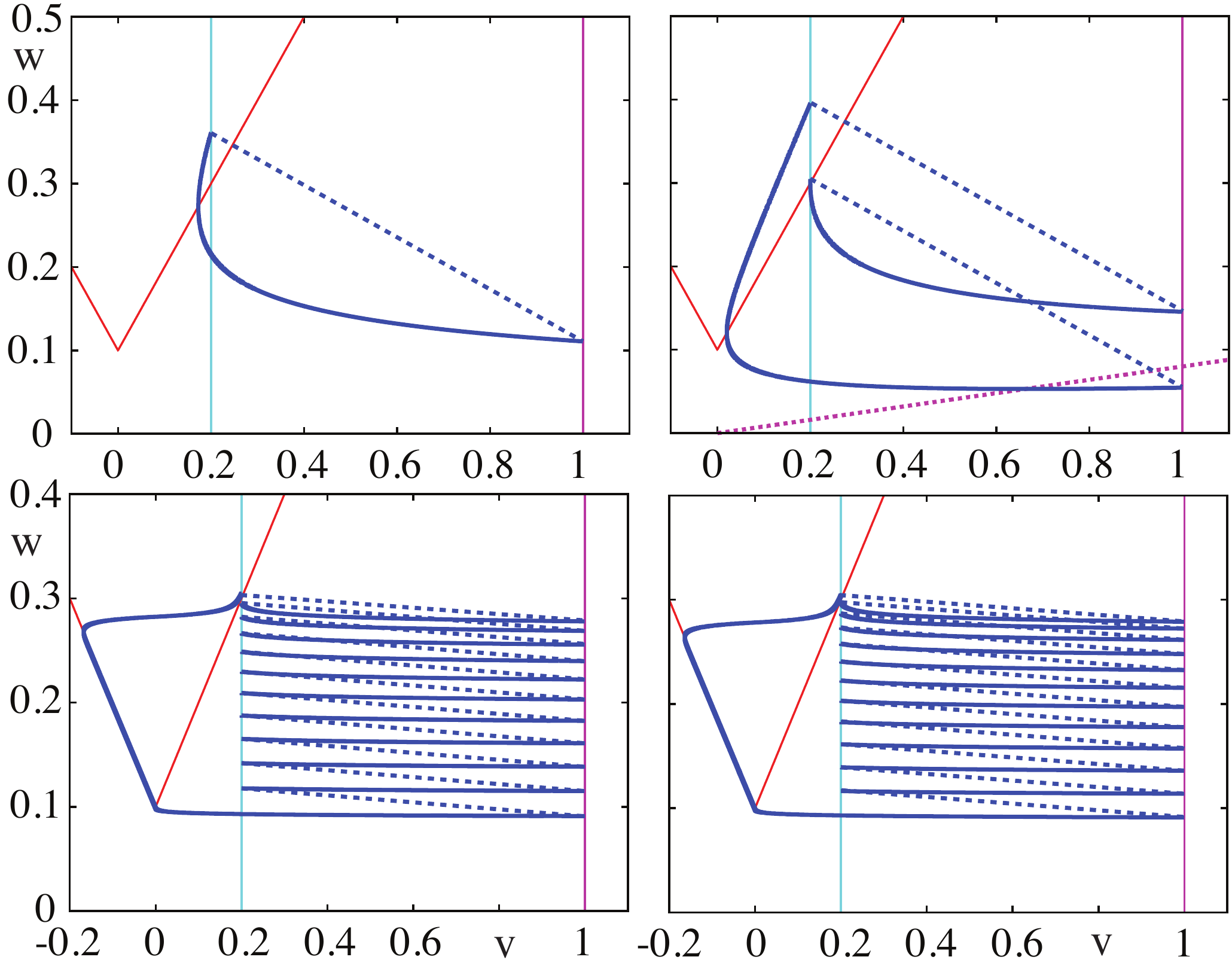}
\caption{Phase plane of the PWL-IF model with $v_\text{th}=1$, $v_\text{r}=0.2$, $a_w=0$, $b_w=-1$, $a_L=-1$, $a_R=1$, $I=0.1$.  The $v$-nullcline is in red (solid) and the $w$-nullcline is in magenta (dashed).
Top Left:  A tonic firing pattern for $\kappa=0.75$ and $\tau=3$.
Top Right:  A spike doublet found for the parameters of tonic firing apart from $a_w=0.08$.
Bottom Left:  A burst firing pattern with $M=11$ spikes per burst for $\kappa=2$ and $\tau=75$.
Bottom Right:  A burst firing pattern with $M=12$ spikes per burst for the parameters of the $11$-spike burst pattern with $\kappa=1.9$
\label{Fig:pwlIF}}
\end{center}
\end{figure}

For a planar PWL IF model that caricatures the adaptive exponential IF model we set $m=2$, $\vec{z}=(v,w)$, $v\leq v_{\text{th}}$, and
\begin{equation}
A_L=\begin{bmatrix}
a_L & -1\\
a_w/\tau & b_w/\tau
\end{bmatrix}, \qquad
A_R=\begin{bmatrix}
a_R & -1\\
a_w/\tau & b_w/\tau
\end{bmatrix},
\nonumber
\end{equation}
with $\vec{c}_L=(I,0)^\mathsf{T}=\vec{c}_R$ for some constant drive $I$.  The switching manifold is prescribed by the choice $v=0$.
Whenever the voltage variable $v$ reaches the firing threshold $v_{\text{th}}$ then the system is reset according to
$\vec{z} \rightarrow \vec{g}(\vec{z}) = (v_{\text{r}}, w +\kappa/\tau)$.  This gives rise to another type of switching manifold that we shall refer to as the \textit{firing manifold}. If we introduce the function $h(\vec{z};a)=v-a$, then the indicator function for \textit{switching} is given by $h(\vec{z};0)$, whilst that for \textit{firing} by $h(\vec{z};v_\text{th})$.
From now on we shall simply refer to the model described above as the \textit{PWL-IF} model.  We note that the vector field defining the PWL-IF model is continuous at $v=0$ (so that there is no need to invoke the Filippov convex method).
A set of phase planes of the model (with nullclines and typical trajectories) for different parameter choices is shown in Fig.~\ref{Fig:pwlIF}.  The PWL-IF model is able to support a tonic firing pattern, which itself is unstable to a period doubling bifurcation as seen in the top part of Fig.~\ref{Fig:pwlIF}.  It can also support a burst firing pattern, which is unstable to \textit{spike-adding} bifurcations as shown in the bottom part of Fig.~\ref{Fig:pwlIF}.
%%%%%%%%%%%%%%%%%%%%%%%%%%%%%%%%%%%%%%%%%%%%%%%%%%%%%%%%%%%%%%%%%%%%%%%%%%%%%%%%%%%%%%%%%%%%%%%%%%%%%%%%%%%%%%%%%%
%%%%%%%%%%%%%%%%%%%%%%%%%%%%%%%%%%%%%%%%%%%%%%%%%%%%%%%%%%%%%%%%%%%%%%%%%%%%%%%%%%%%%%%%%%%%%%%%%%%%%%%%%%%%%%%%%%
Each of the periodic orbits shown in Fig.~\ref{Fig:pwlIF} can be obtained without recourse to the numerical evolution of the underlying nonsmooth flow for the PWL-IF model.  Rather the shape of the orbits can be determined explicitly by exploiting the linearity of the model between switching and firing manifolds.  In any such region the model is given by a linear system of ordinary differential equations (ODEs) of the form $\dot{z}=A \vec{z}+ \vec{c}$ for appropriate choices of $A$ and $\vec{c}$.  The explicit solution can be obtained in terms of initial data using matrix exponentials as
\begin{equation}
%\vec{z}(t) = G(t) \vec{z}(0) + J(t) \vec{c},
\vec{z}(t) = \e^{At} \vec{z}(0) + A^{-1}[\e^{At}-I_m] \vec{c} .
\label{solu}
\end{equation}
%where $G(t)= \e^{At}$, $J(t) = A^{-1}[G(t)-I_m]$, and $I_m$ is the $m \times m $ identity matrix.

For example, to determine the tonic firing pattern shown in Fig.~\ref{Fig:pwlIF}, which only ever visits the region of phase-space with $\vec{z} \in S_R$ we would simply use the formula for (\ref{solu}) with $A=A_R$ and $\vec{c}=\vec{c}_R$ to determine the \textit{time-of-flight} $\Delta$ according to $v(\Delta)=v_\text{th}$ with $(v(0),w(0))=(v_r,w_0)$, with $w_0$ determined self-consistently according to $w_0=w(\Delta)+ \kappa/\tau$.
The resulting pair of nonlinear equations for $(\Delta,w_0)$ can in general be solved with a numerical root finding scheme (and choosing the solution with  the largest value of $w_0$).  Thus although we can eliminate the need to numerically solve an ODE, there is still some need for root-finding.  To determine the burst patterns shown in Fig.~\ref{Fig:pwlIF}, with $M$ spikes per burst, would require the simultaneous solution of $M+2$ equations (one for the value of $w$ on the switching manifold, and the others to determine the times-of-flight between consecutive firing events, remembering that for this orbit the trajectory also visits the region $S_L$).  However, once the orbit is determined in this way the Floquet theory for stability simplifies considerably.  To see how this simplification arises it is enough to focus on the stability of simple tonic firing pattern.  In this case the Jacobian of the orbit is the constant matrix $A_R$, and is independent of the orbit shape.  Thus small perturbations to the periodic orbit $\delta \vec{z}(t)$ can be simply constructed according to $\delta \vec{z}(t)  = \exp(A_R t) \delta \vec{z}(0)$.  The caveat being that this result only holds away from switching or firing events.  To propagate perturbations properly through switching and firing manifolds we make use of a \textit{saltation} matrix \cite{Aizerman1958}.  A derivation of the saltation matrix in the context of this paper is given in Appendix A.  Thus for our nonsmooth system the evolution of perturbations to the orbit over one period is given by $\delta \vec{z}(\Delta)  = \Psi \delta \vec{z}(0)$, where $\Psi=K(\Delta) \exp(A_R \Delta)$ and $K$ is the saltation matrix:
\begin{equation}
K(t) = \begin{bmatrix}
\dot{v}(t^+)/\dot{v}(t^-) & 0 \\
(\dot{w}(t^+)-\dot{w}(t^-))/\dot{v}(t^-) & 1
\end{bmatrix} .
\nonumber
\end{equation}
The orbit will be stable provided that the eigenvalues of $\Psi$ lie within the unit disc.  Since one of these will be unity (reflecting time-translation invariance) we have that the orbit is stable provided $\text{Re}\,r <0$, where the \textit{Floquet exponent} $r$ can be calculated explicitly as
\begin{align}
r &= \operatorname{Tr} A_R + \frac{1}{\Delta} \ln \frac{\dot{v}(\Delta^+)}{\dot{v}(\Delta^-)} \nonumber \\
&=a_R+\frac{b_w}{\tau} + \frac{1}{\Delta} \ln \frac{a_R v_\text{r}-w(0)+I}{a_R v_\text{th}-w(\Delta)+I} . \nonumber
\end{align}
Here we have used the fact that $\det \Psi = 1 \times \lambda$, and defined the non-trivial eigenvalue $\lambda$ of $\Psi$ as $\lambda = \exp(r \Delta)$.

\begin{figure}
\begin{center}
\includegraphics[width=6cm]{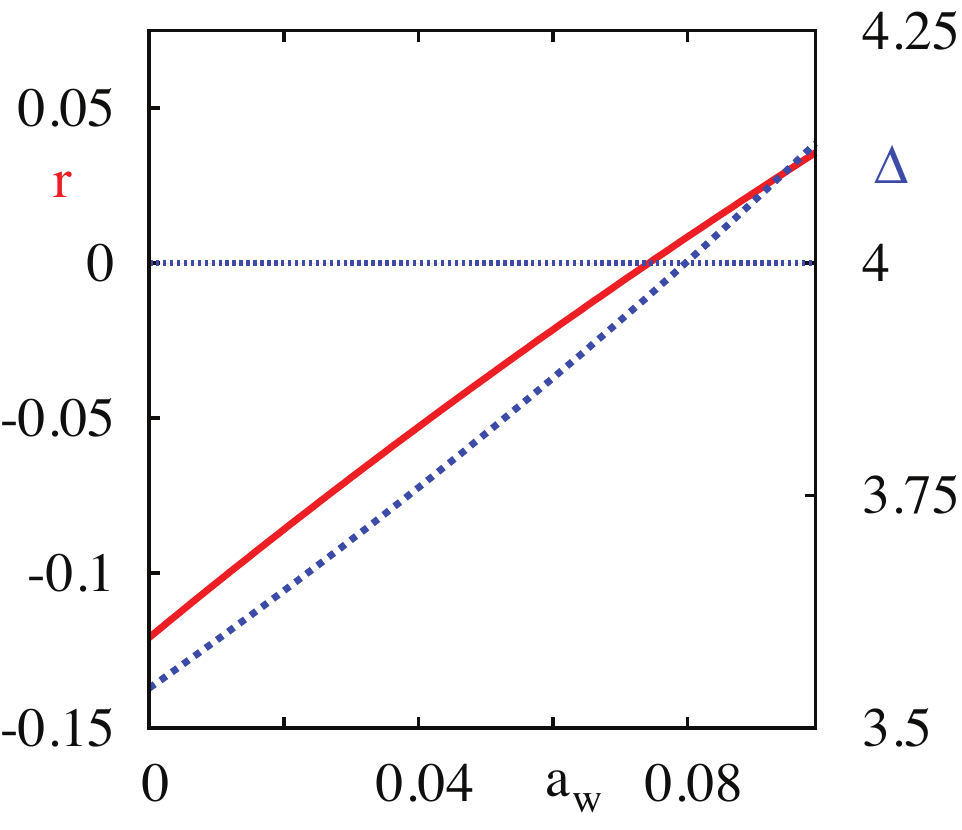}
\caption{A plot of the Floquet exponent $r$ (red solid) for the tonic orbit shown in Fig.~\ref{Fig:pwlIF} as a function of $a_w$.  Here we see that a bifurcation, defined by $\text{Re} \, r=0$, occurs as $a_w$ increases through $0.075$. The corresponding value of $\text{Im}\, r$ is found to be $\pi/\Delta$, signalling a period doubling instability (confirmed in direct numerical simulations).  Also shown is the variation of the period $\Delta$ (blue dashed) using the right hand vertical axis.
\label{Fig:sigma}}
\end{center}
\end{figure}

In Fig.~\ref{Fig:sigma} we show a plot of $\text{Re} \, r$ as a function of the parameter $a_w$.  This predicts the point of a period doubling instability, and is found to be in excellent agreement with direct numerical simulations.  The Floquet exponent for more complicated orbits (that burst and visit both $S_L$ and $S_R$) can also be easily constructed using a revised structure for $\Psi$.  For example for a burst pattern with $M$ spikes of the type shown in Fig.~\ref{Fig:pwlIF} then we would have that
\begin{align}
\Psi &=  K(T_{M+2})\e^{A_L \Delta_{M+2}}K(T_{M+1})\e^{A_R \Delta_{M+1}} \nonumber \\
& \times K(T_M)\e^{A_R \Delta_M} \ldots K(T_2)\e^{A_R \Delta_2}K(T_1)\e^{A_R \Delta_1} , \nonumber
\end{align}
where $\Delta_i$ indicate various times-of-flight, and the $T_i$ the times at which firing or switching events occur.
%%%%%%%%%%%%%%%%%%%%%%%%%%%%%%%%%%%%%%%%%%%%%%%%%%%%%%%%%%%%%%%%%%%%%%%%%%%%%%%%%%%%%%%%%%%%%%%%%%%%%%%%%%%%%%%%%%
%%%%%%%%%%%%%%%%%%%%%%%%%%%%%%%%%%%%%%%%%%%%%%%%%%%%%%%%%%%%%%%%%%%%%%%%%%%%%%%%%%%%%%%%%%%%%%%%%%%%%%%%%%%%%%%%%%

We now consider a network of $N$ synaptically coupled PWL-IF neurons with the introduction of an index $i=1,\ldots,N$ (labelling each node) and the replacement of the constant external drive by a time-dependent forcing such that
$I \rightarrow I + \sigma \sum_j W_{ij} s_j(t)$.  Here the synaptic input from neuron $j$ takes the standard \textit{event-driven} form
\begin{equation}
s_j(t) =  \sum_{p \in \Zset} \eta(t-T_j^p),
\nonumber
\end{equation}
where $T_j^p$ denotes the $p$th firing time of neuron $j$ and $\eta$ describes the shape of the post-synaptic response.  Other forms of \textit{state-dependent} synaptic coupling are also possible, such as the fast threshold modulation type used by Belykh and Hasler to study clusters in networks of bursting neurons \cite{Belykh2011}, though would require a further PWL reduction before they could be studied with the method presented here.
For a balanced network the existence of the synchronous network state is independent of any of the parameters describing the synaptic interaction.  This means that the form of synaptic coupling cannot induce any nonsmooth bifurcations, such as grazing, which occur if an orbit tangentially touches the firing threshold.
Here we shall focus on the common choice of a continuous $\alpha$-function, so that $\eta(t) = \alpha^2 t \e^{-\alpha t} \Theta(t)$, where $\Theta$ is a Heaviside step function.  In this case we may also write $s_i(t)$ as the solution to an impulsively forced linear system:
\begin{align}
\left ( 1 + \frac{1}{\alpha} \FD{}{t} \right ) s_i &= u_i, \nonumber\\
\left ( 1 + \frac{1}{\alpha} \FD{}{t} \right ) u_i &= \sum_{p \in \Zset} \delta(t-T_i^p).
\label{history}
\end{align}
Exploiting the linearity of the synaptic dynamics between firing events we may succinctly write the network model in the form $\dot{\vec{z}}_i = \vec{F}(\vec{z}_i)$, where $\vec{z}_i = (v_i,w_i,s_i,u_i)$, and $\vec{F}$ has the form of (\ref{zdot}) with
\begin{equation}
A_{L,R}=\begin{bmatrix}
a_{L,R} & -1 & 0 & 0 \\
a_w/\tau & b_w/\tau & 0 & 0 \\
0 & 0 & -\alpha & \alpha \\
0 & 0 & 0 & -\alpha
\end{bmatrix},
\nonumber
\end{equation}
and $\vec{c}_L=(I,0,0,0)^\mathsf{T}=\vec{c}_R$, with $\vec{z}_i \rightarrow \vec{g}(\vec{z}_i) = (v_{\text{r}}, w_i +\kappa/\tau, s_i, u_i+ \alpha)$ whenever $h(\vec{z}_i;v_\text{th})=0$.  The vector function that specifies the interaction is given by $\vec{H}(\vec{z}_i) = (s_i,0,0,0)^\mathsf{T}$.
The MSF approach for a smooth network gives rise to a Floquet problem with a Jacobian $D \vec{F}(\vec{s}) + \sigma \lambda_l D \vec{H}(\vec{s})$, where $\lambda_l$ is an eigenvalue of the coupling matrix $W$.  The corresponding Jacobian for the PWL-IF network is $A_{L,R} + \sigma \lambda_l D \vec{H}$, with the label $L$ or $R$ chosen according to whether the synchronous orbit is in $S_L$ or $S_R$.  Moreover, $D \vec{H}$ is now a constant matrix with $[D \vec{H}]_{ij} = 1$ if $i=1$ and $j=3$, and is $0$ otherwise.

\begin{figure}
\begin{center}
\includegraphics[width=8cm]{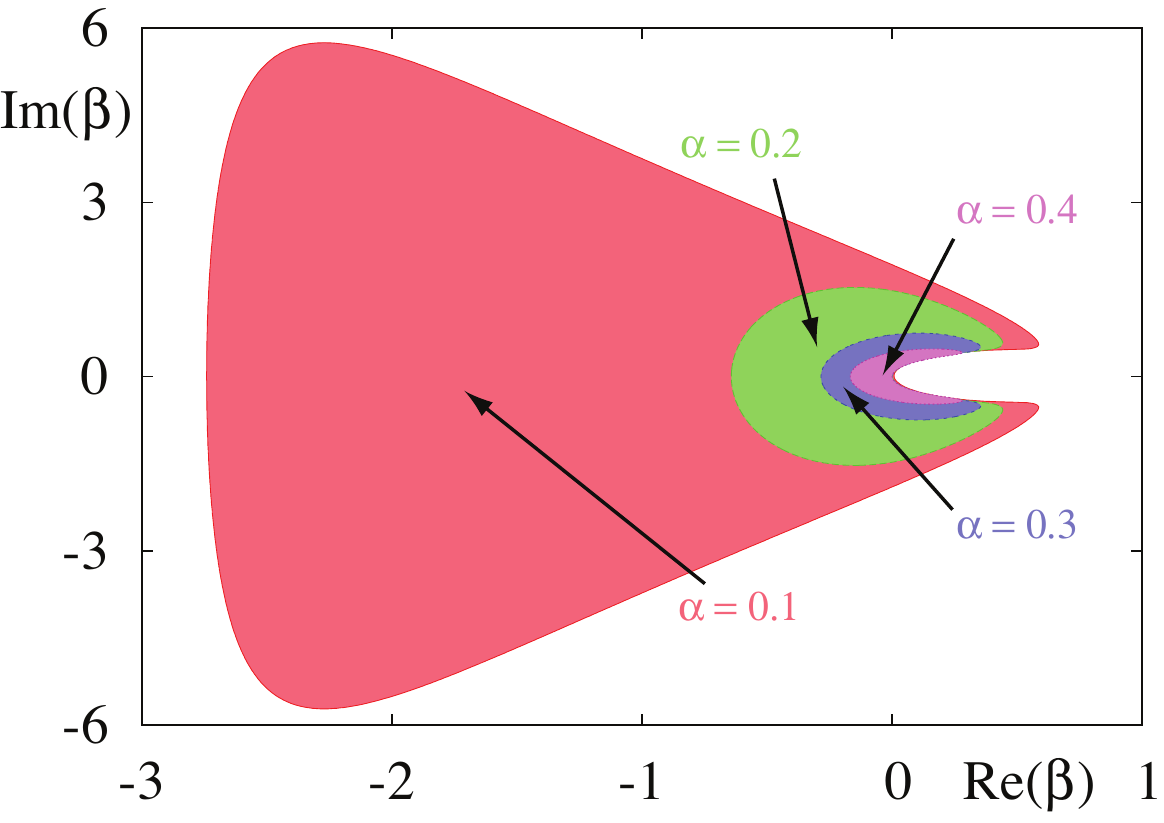}
\caption{The MSF for a PWL-IF network with a synchronous tonic orbit.  Parameters as in Fig.~\ref{Fig:pwlIF} top left panel.  The shaded regions indicate where $\text{MSF}<0$ for various values of the synaptic rate parameter $\alpha$.  The largest region is for $\alpha=0.1$, with correspondingly smaller areas for $\alpha=0.2,0.3,0.4$ respectively.  The synchronous solution is stable provided all of the eigenvalues of $\sigma W$ lie within a shaded area for a given value of $\alpha$.
\label{Fig:MSF}}
\end{center}
\end{figure}

The propagation of (linearised) trajectories through a switch or a firing event is achieved with the use of a saltation matrix $K\in \Rset^{m \times m}$, and we write this as $\vec{U}^+ = (I_N \otimes K) \vec{U}^-$, where $\vec{U}^-$ and $\vec{U}^+$ are respectively the system states just before and just after the saltation.  In transformed coordinates (and see equation (\ref{Vdot})) it is simple to establish that $\vec{V}^+ = (I_N \otimes K) \vec{V}^-$, so that saltation also acts blockwise with $\xi_l^+ = K \xi_l^-$.
Thus the synchronous solution is stable (following the same arguments as for a single node) if all of the eigenvalues of a set of matrices $\Psi_l$, $l=1,\ldots,N$ lie within the unit disc, excluding the one that arises from time-translation invariance (with a value $+1$).
For example for a synchronous tonic orbit of the type discussed above,
\begin{equation}
\Psi_l = K(\Delta) \exp \{ (A_R + \sigma \lambda_l D \vec{H}) \Delta \} .
\nonumber
\end{equation}
Here the saltation matrix is given by (and see Appendix A)
\begin{equation}
K(t) = \begin{bmatrix}
\dot{v}(t^+)/\dot{v}(t^-) & 0 & 0 & 0\\
(\dot{w}(t^+)-\dot{w}(t^-))/\dot{v}(t^-) & 1 & 0 & 0 \\
(\dot{s}(t^+)-\dot{s}(t^-))/\dot{v}(t^-) & 0 & 1 & 0\\
(\dot{u}(t^+)-\dot{u}(t^-))/\dot{v}(t^-) & 0 & 0 & 1\\
\end{bmatrix} .
\nonumber
\end{equation}
The periodic trajectory for $\vec{z}(t)$ is subject to the constraints $v(\Delta)=v_\text{th}$, $w(0)=w(\Delta)+\kappa/\tau$, $s(0)=s(\Delta)$, and $u(0)=u(\Delta)+\alpha$.
If we denote an eigenvalue of $K(\Delta) \exp \{ (A_R + \beta D \vec{H}) \Delta \}$, $\beta \in \Cset$, by $\gamma(\beta)$ then the MSF is the largest number in the set $\text{Re} \, ( \ln \gamma(\beta) )/\Delta$, and the synchronous state is stable if the MSF is negative at all the points where $\beta=\sigma \lambda_l$.
In Fig.~\ref{Fig:MSF} we show a plot of the MSF for the synchronous tonic orbit.
This solution will be stable provided all of the eigenvalues of $\sigma W$ lie within the shaded area  shown in Fig.~\ref{Fig:MSF} (for a given $\alpha$). For a positive semi-definite connectivity matrix then we see from Fig.~\ref{Fig:MSF} that the synchronous solution is unstable for $\sigma >0$.  However, for $\sigma<0$ the same network can support a stable synchronous orbit for some sufficiently small value of $|\sigma|$.
%%%%%%%%%%%%%%%%%%%%%%%%%%%%%%%%%%%%%%%%%%%%%%%%%%%%%%%%%%%%%%%%%%%%%%%%%%%%%%%%%%%%%%%%%%%%%%%%%%%%%%%%%%%%%%%%%%
%%%%%%%%%%%%%%%%%%%%%%%%%%%%%%%%%%%%%%%%%%%%%%%%%%%%%%%%%%%%%%%%%%%%%%%%%%%%%%%%%%%%%%%%%%%%%%%%%%%%%%%%%%%%%%%%%%

We note here that for the choice of an exponential synapse as originally considered in \cite{Ladenbauer2013}, with $\eta(t) = \alpha \e^{-\alpha t} \Theta(t)$, then synaptic interactions are \textit{discontinuous} and the order in which perturbed components of the state vector cross the firing threshold becomes important \cite{Coombes2000}.  In this case the approach above, valid only for continuous interactions, must be modified (although this was not considered in \cite{Ladenbauer2013}).  We shall treat this mathematically interesting case in another paper.

%%%%%%%%%%%%%%%%%%%%%%%%%%%%%%%%%%%%%%%%%%%%%%%%%%%%%%%%%%%%%%%%%%%%%%%%%%%%%%%%%%%%%%%%%%%%%%%%%%%%%%%%%%%%%%%%%%
\begin{figure*}[hbtp]
\begin{center}
\includegraphics[width=16cm]{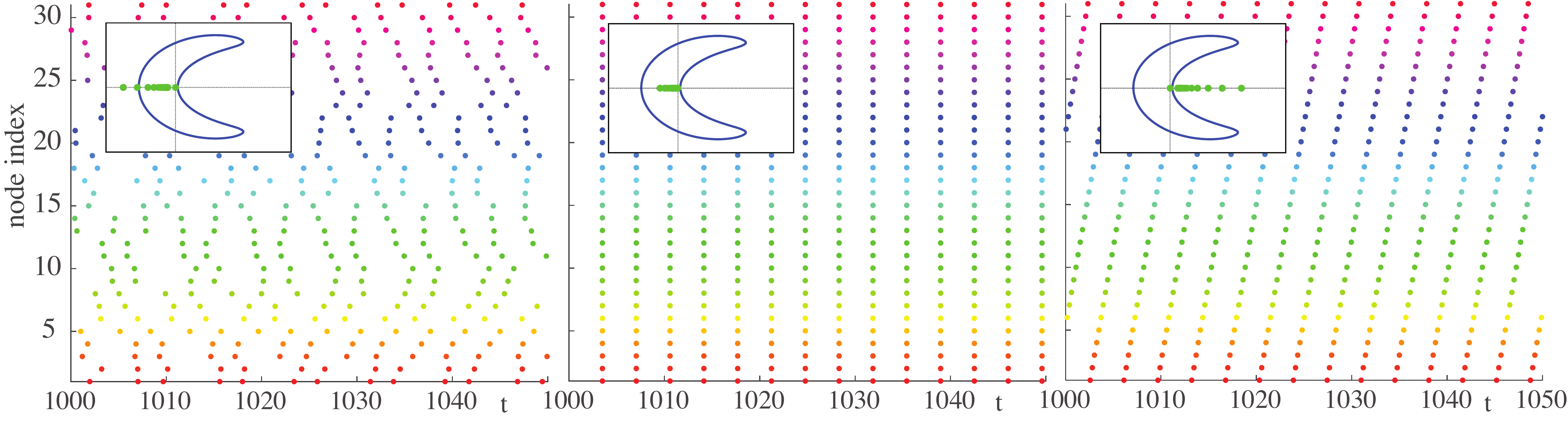}
\caption{Raster plot of spike times from direct numerical simulations of a PWL-IF network with $N=31$, $d=3$, and $\alpha=0.4$.  The insets show a plot of the MSF with the eigenvalues of $\sigma W$ superimposed.
Left:  $\sigma=-0.1$, and synchrony is unstable.
Middle: $\sigma=-0.025$, and synchrony is stable.
Right: $\sigma=0.1$, and synchrony is unstable.
The predicted instability borders (at $\sigma=0$ and $\sigma \simeq -0.05$) are in good agreement with the predictions from the MSF analysis.  The typical pattern of firing activity beyond an instability of the synchronous state for $\sigma >0$ is found to be a periodic travelling wave, whilst for $\sigma <0$ a spatio-temporal pattern emerges via a period doubling instability of the firing times.
\label{Fig:pwlIFsim}}
\end{center}
\end{figure*}
%%%%%%%%%%%%%%%%%%%%%%%%%%%%%%%%%%%%%%%%%%%%%%%%%%%%%%%%%%%%%%%%%%%%%%%%%%%%%%%%%%%%%%%%%%%%%%%%%%%%%%%%%%%%%%%%%%

As a particular realisation of a network architecture that guarantees synchrony we choose a balanced ring network with $N$ odd and $W_{ij} = W(|i-j|)$, with distances calculated modulo $(N-1)/2$, and $w(x) = (1-a|x|/d) \e^{-|x|/d}$.  Here the parameter $a$ is chosen such that, for a given value of $N$ and a scale $d$, $\sum_{j=1}^N W_{ij}=0$.  The circulant structure of this matrix means that it has normalised eigenvectors given by $\vec{e}_l = (1,\omega_l, \omega_l^2, \ldots, w_l^{N-1})/\sqrt{N}$, where $l=0,\ldots, N-1$, and $\omega_l=\exp(2 \pi i l /N)$ are the $N$th roots of unity.  The eigenvalues of the (symmetric) connectivity matrix are real and given by
\begin{equation}
\lambda_l = \sum_{j=0}^{N-1} W(|j|) \omega_l^j.
\nonumber
\end{equation}
We note that the balance condition enforces $\lambda_0=0$.  We also have that $\lambda_{N-l}=\lambda_l$ for $l=1,\ldots, (N-1)/2$, so that any excited pattern is given by a combination $\vec{e}_m+\vec{e}_{-m} = 2 \text{Re} \, \vec{e}_m$ for some $1\leq m \leq (N-1)/2$.  Given the shape of the MSF function shown in Fig.~\ref{Fig:MSF} the value of $m$ is determined by $\lambda_m= \max_{l} \lambda_l$.  In Fig.~\ref{Fig:pwlIFsim} we compare simulations of a network against the predictions of the MSF.  When the network eigenvalues lie within the region where the MSF is negative, then small perturbations to synchronous initial data die away and the system settles to a synchronous periodic orbit as expected.  When one of the eigenvalues crosses the zero level set of the MSF (from negative to positive) we see two different types of instability emerge.  One is the emergence of a spatio-temporal pattern of spike-doublets, arising because an eigenvalue of $\Psi_l$ leaves the unit disc at $-1$ (a period-doubling bifurcation), and another giving rise to periodic travelling wave (with asynchronous firing) because an eigenvalue of $\Psi_l$ leaves the unit disc at $+1$ (tangent bifurcation).

\section{\label{sec:Cluster}Cluster states and network structure}

We now move on to discuss the more general phenomenon of cluster synchronisation, where different groups of oscillators are exactly synchronised but there is no exact synchrony between the groups. We focus on networks with symmetry, where cluster states arise very naturally (although clustering can also occur in networks without symmetry where synchronous nodes have synchronous input patterns \cite{Golubitsky2016}). For networks of identical oscillators with dynamics described by \eqref{network1}, a symmetry of the network is a permutation $\gamma$ of the nodes which leaves the network equations unchanged. That is, if we denote $\vec{x}= (\vec{x}_1, \ldots, \vec{x}_N) \in \mathbb{R}^{Nm}$, $M_\gamma$ the $N\times N$ permutation matrix for the permutation $\gamma$ and \[\vec{G}_i(\vec{x}) = \vec{F}(\vec{x}_i) - \sigma \sum_{j=1}^N \mathcal{G}_{ij} \vec{H}(\vec{x}_j),\] then the network dynamics \eqref{network1} are given by $\dot{\vec{x}} = \vec{G}(\vec{x})$ for $\vec{G}= (\vec{G}_1, \ldots, \vec{G}_N)^\mathsf{T}$ where the vector field $\vec{G}$ satisfies
\begin{equation} \label{eq:equivariance}
\vec{G}((M_\gamma \otimes I_m) \vec{x}) =(M_\gamma \otimes I_m) \vec{G}(\vec{x}).\end{equation}
 This results in the condition that $\gamma$ is a symmetry of the network if $M_\gamma \mathcal{G} = \mathcal{G} M_\gamma$ (or equally $M_\gamma$ commutes with the connectivity matrix). The network symmetries form a group $\Gamma \subseteq S_N$. Many real-world networks have been shown to have a high degree of symmetry, arising from locally tree-like structures produced by natural growth of the network \cite{MacArthur2008}. Determining symmetries of such large and complex networks is impossible by inspection. Indeed even networks with a small number of nodes can have a large symmetry group \cite{Pecora2014}. However, the computations required to determine the symmetry group of a given network are easily implemented using computational algebra routines \cite{GAP4, sagemath}.

Recent work of Pecora \textit{et al}. \cite{Pecora2014} and Sorrentino \textit{et al}. \cite{Sorrentino2016} has demonstrated how the ideas behind the MSF can be combined with group theoretical techniques used in the study of symmetric dynamical systems to analyse the stability of cluster states within symmetric networks of dynamical units. The approach taken in \cite{Pecora2014} and \cite{Sorrentino2016} continues to focus on networks of coupled identical oscillators whose dynamics are given by \eqref{network1}. This approach is equivalent to a restriction (to a particular form of admissible network equations \eqref{network1} in the special case of identical oscillators and one type of bidirectional coupling) of the more general theory of patterns of synchrony in coupled cell systems developed by Golubitsky and Stewart and collaborators and recently reviewed in \cite{Golubitsky2016}. Using terminology from the more general theory, here we review for the particular network dynamics \eqref{network1} how network structure can be used to determine a catalogue of possible cluster states. We also discuss techniques for determining the stability of any given cluster state, utilising methods from computational group theory. The exposition below is inspired by that of \cite{Sorrentino2016} and is presented as a recipe for how to extend the standard MSF approach to treat networks of PWL oscillators.

Many of the possible cluster states which a given network may support arise from the symmetry of the network. We first describe how these cluster states may be determined before highlighting the algorithm of Sorrentino \textit{et al}. \cite{Sorrentino2016} which can be used to determine additional cluster states resulting from the particular choice of Laplacian coupling.

\subsection{\label{sec:ClustersSymmetry}Cluster states from network symmetries}

Consider a network of $N$ identical oscillators with one type of bidirectional coupling whose dynamics are given by \eqref{network1} and which has symmetry group $\Gamma$.  As a consequence of the dynamics satisfying the equivariance condition \eqref{eq:equivariance}, if $\vec{x}(t)$ is a solution of \eqref{network1} (equilibrium or periodic state) then for any permutation $\gamma \in \Gamma$,
\[ \FD{}{t} \gamma \vec{x} = \gamma \FD{\vec{x}}{t} = \gamma \vec{G}(\vec{x}) = \vec{G}(\gamma \vec{x}),\] so $\gamma \vec{x}(t)$ is also a solution. Here $\gamma$ acts on $\vec{x}$ via the permutation matrix $(M_\gamma \times I_m)$. Thus solutions occur as group orbits, $\{ \gamma \vec{x} \ : \ \gamma \in \Gamma\}$. The isotropy subgroup of a solution is the subgroup $\Sigma \subseteq \Gamma$ given by
\begin{equation} \nonumber \Sigma = \{ \gamma \in \Gamma \ : \ \gamma \vec{x}(t) = \vec{x}(t), \ \forall t \in \mathbb{R}\}.\end{equation} This is the group of (spatial) symmetries of the solution $\vec{x}(t)$ and the largest subgroup of $\Gamma$ under which the solution is invariant. Solutions which lie on the same group orbit have conjugate isotropy subgroups and the same existence and stability properties. Given a subgroup $\Sigma \subseteq \Gamma$ we can define its fixed-point subspace
\[ \mbox{Fix}(\Sigma) = \{ \vec{x} \ : \ \gamma \vec{x} = \vec{x},\  \forall \gamma \in \Sigma\}.\]
Fixed-point subspaces are flow invariant: Let $\gamma \in \Sigma$ and $\vec{x} \in \mbox{Fix}(\Sigma)$. Then
\[ \gamma \vec{G}(\vec{x}) = \vec{G}(\gamma \vec{x}) = \vec{G}(\vec{x}),\] so $\vec{G}(\vec{x}) \in \mbox{Fix}(\Sigma)$.

Suppose that %$\Sigma$ is the isotropy subgroup of some solution $\vec{x} = (\vec{x}_1, \ldots, \vec{x}_N)$.
 $\Sigma$ is any subgroup of $\Gamma$. The orbit under $\Sigma$ of the node $i$ is the set $\{ \gamma(i) \ : \ \gamma \in \Sigma\}$. The orbits permute subsets of nodes among each other and in this way partition the nodes into clusters. Nodes which lie on the same orbit (in the same cluster) have synchronised dynamics, $\vec{x}_{\gamma(i)} \equiv \vec{x}_i$ for any $\gamma \in \Sigma$ - see \cite[Thm III.2]{Golubitsky2016}. Also, the synchronised state for each cluster is flow-invariant. Thus, to enumerate all possible synchronised cluster states of a given network which are due to network symmetries, we need to determine the network symmetry group $\Gamma$ and all of its subgroups. There will be one type of cluster state for each isotropy subgroup for the action of $\Gamma$ on the node space up to conjugacy.

\medskip
\noindent{\textbf{Example: A five-node network}}

\begin{figure}
\begin{center}
\includegraphics[width=2.5cm]{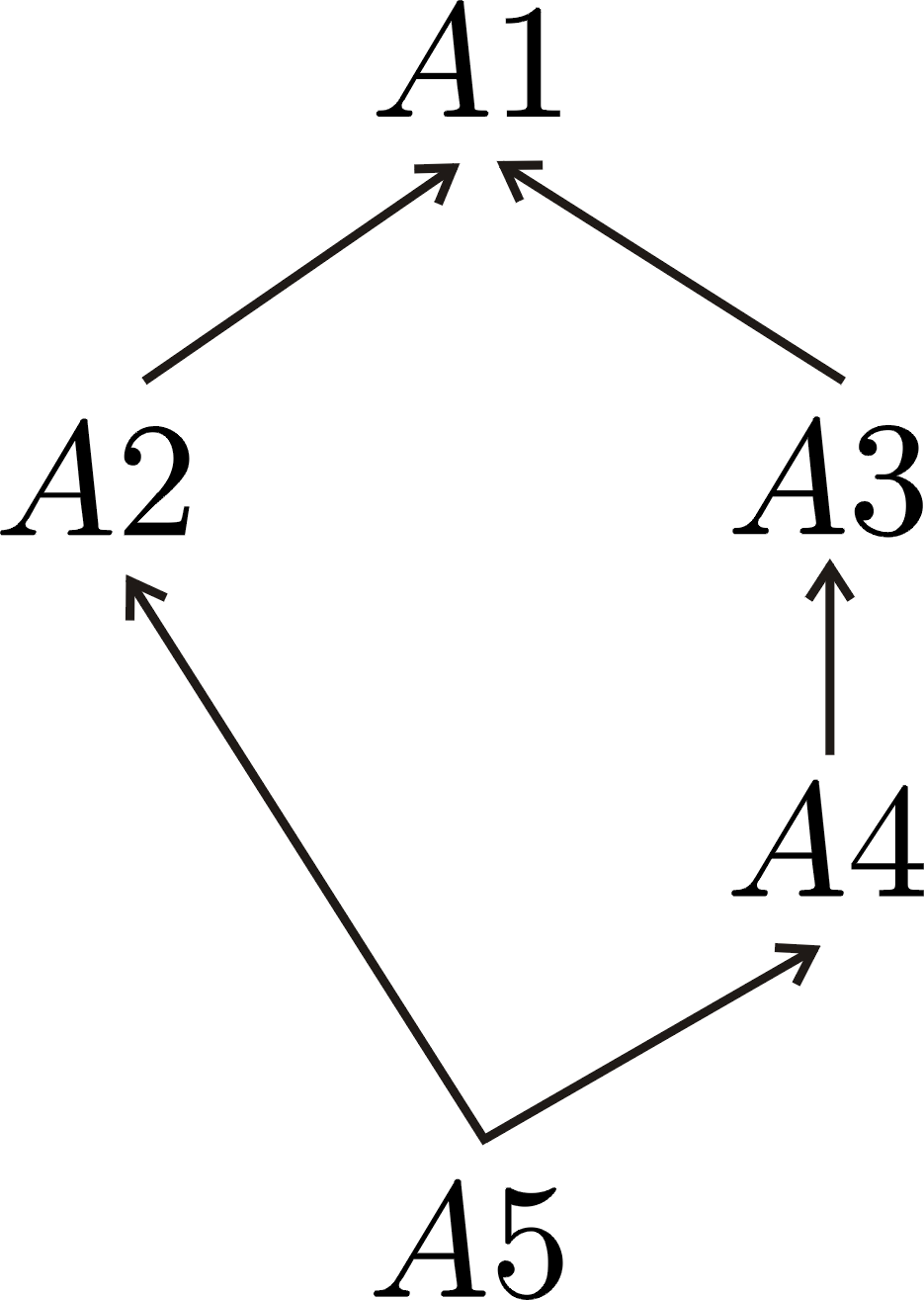}
\caption{The lattice of cluster states given by network symmetry for the five-node network with connectivity \eqref{eq:Laplacian5}. Arrows indicate inclusion.
\label{Fig:lattice}}
\end{center}
\end{figure}

\smallskip
\noindent
Consider the example five-node network studied by Sorrentino \textit{et al}. \cite{Sorrentino2016} which has graph Laplacian matrix
 \begin{equation}\label{eq:Laplacian5}
 \mathcal{G} = \begin{pmatrix*}[r] 3 &-1 &0 & -1 & -1\\ -1 & 3 &-1 & 0 & -1  \\ 0 &-1 &3 & -1 & -1\\-1 & 0 &-1 & 3 & -1  \\ -1 & -1 & -1 & -1 & 4\end{pmatrix*}.
 \end{equation}
  This network has symmetry group $\Gamma \cong \mathbf{D}_4$ generated by the permutations
 \[ \rho= (1234), \quad \pi=(24).\] The conjugacy classes of isotropy subgroups for this network and corresponding cluster states are given in Table \ref{tab:clusters5}. Nodes belonging to the same cluster are the same colour. These are all possible cluster states arising from network symmetries. Figure \ref{Fig:lattice} shows the lattice of isotropy subgroups (cluster states) where arrows indicate (up to conjugacy) inclusion of the symmetries of the cluster state at the tail of the arrow within the symmetry group of the cluster state at the head of the arrow. This gives an indication of likely symmetry-breaking bifurcations between these cluster states.  \hfill $\Diamond$
\medskip

\begin{table*}[hbtp]
\begin{center}
\begin{tabular}{|c|c|c|c|}
\hline
Label & Isotropy & Cluster state & $\mathcal{G}^\prime$\\
\hline
A1 & $\Gamma \cong \mathbf{D}_4$ & \raisebox{-.5\height}{\includegraphics[height=2.3cm]{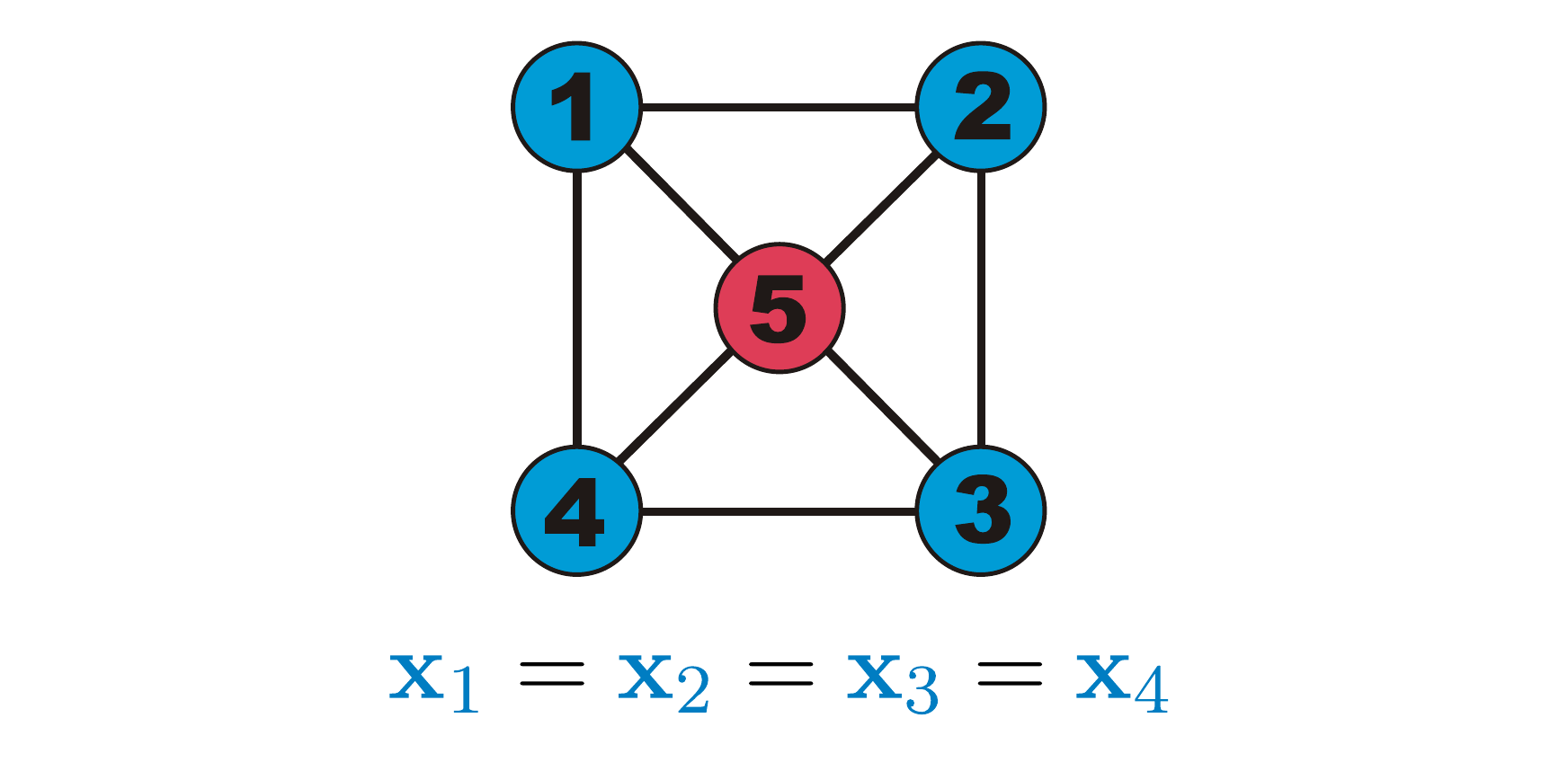}} & $\begin{pmatrix*}[r] 4 &2 &0 & 0 & 0\\ 2 & 1 &0 & 0 & 0  \\ 0 &0 &5 & 0 & 0\\0 & 0 &0 & 3 & 0  \\ 0 & 0 & 0 & 0 & 3\end{pmatrix*}$\\
\hline
A2 & $\mathbb{Z}_2(\pi\rho) \cong \mathbb{Z}_2(\pi\rho^3)$ &  \raisebox{-.5\height}{\includegraphics[height=2.3cm]{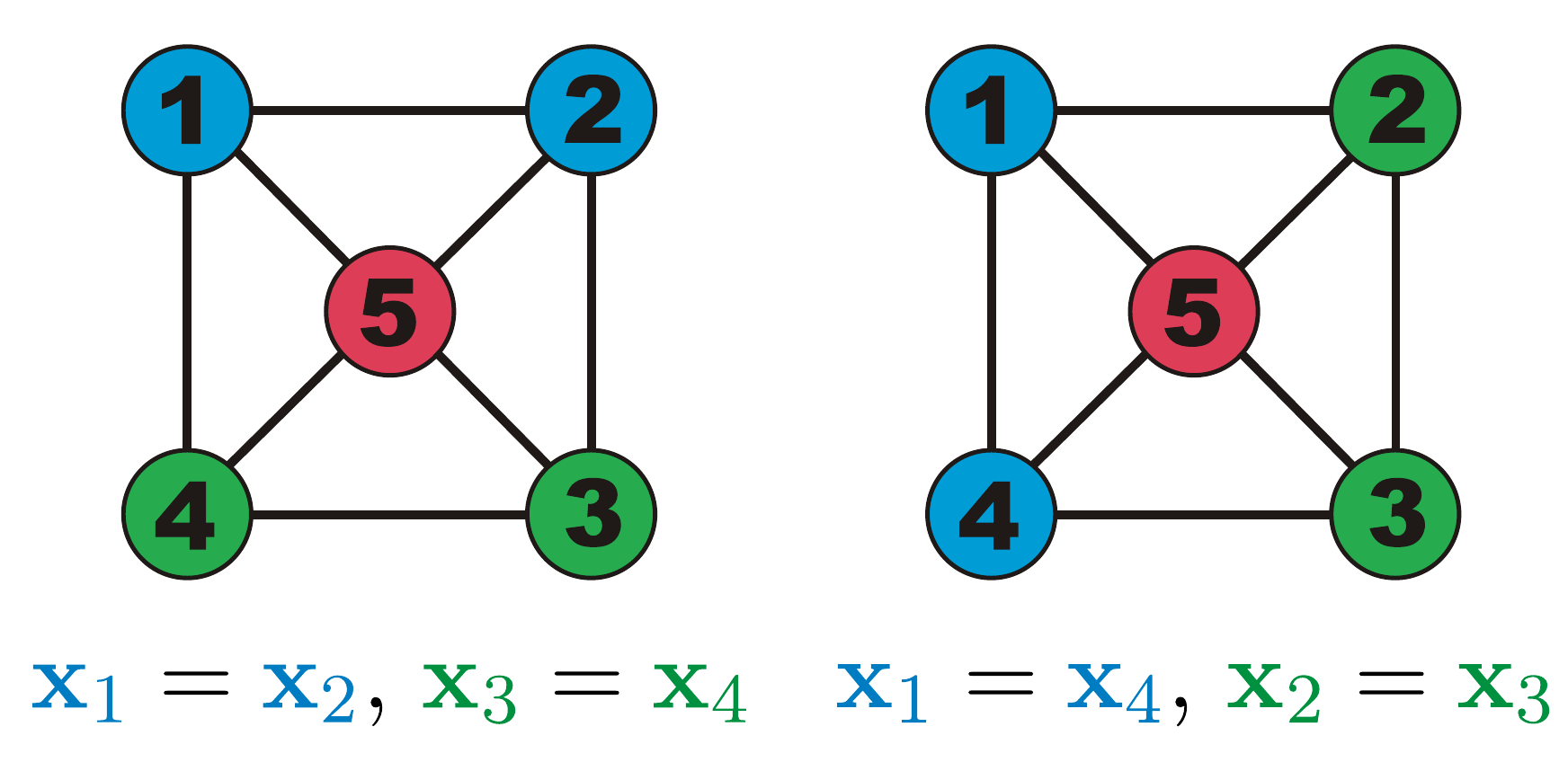}} & $\begin{pmatrix*}[r] 4 &\sqrt{2} &-\sqrt{2} & 0 & 0\\ \sqrt{2} & 2 &1 & 0 & 0  \\ -\sqrt{2} &1 &2 & 0 & 0\\0 & 0 &0 & 4 & 1  \\ 0 & 0 & 0 & 1 & 4\end{pmatrix*}$\\
\hline
A3 & $\mathbf{D}_2(\pi, \pi\rho^2)$ & \raisebox{-.5\height}{\includegraphics[height=2.3cm]{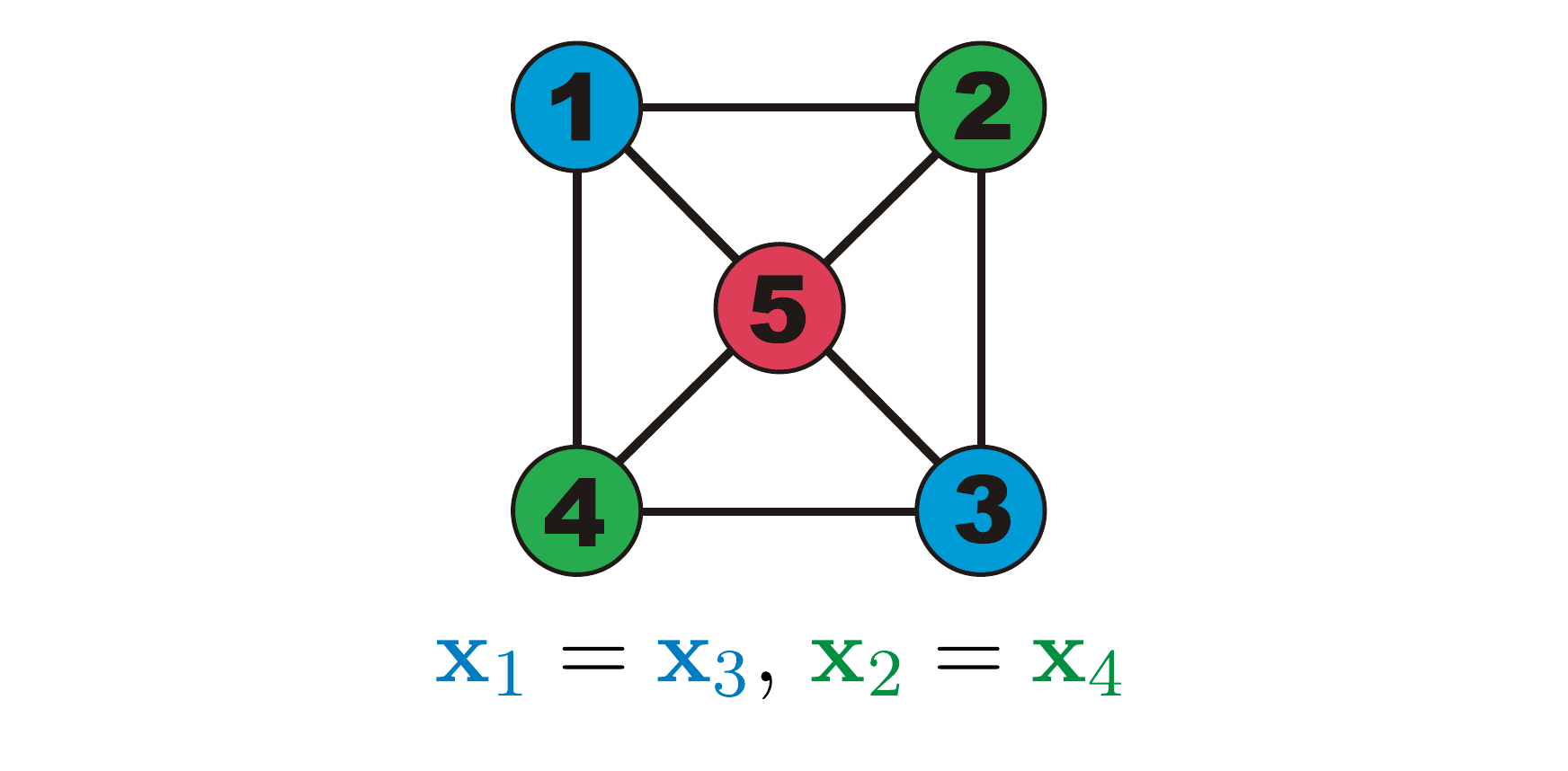}} & $\begin{pmatrix*}[r] 4 &\sqrt{2} &-\sqrt{2} & 0 & 0\\ \sqrt{2} & 3 &2 & 0 & 0  \\ -\sqrt{2} &2 &3 & 0 & 0\\0 & 0 &0 & 3 & 0  \\ 0 & 0 & 0 & 0 & 3\end{pmatrix*}$\\
\hline
A4 & $\mathbb{Z}_2(\pi) \cong \mathbb{Z}_2(\pi \rho^2)$ & \raisebox{-.5\height}{\includegraphics[height=2.3cm]{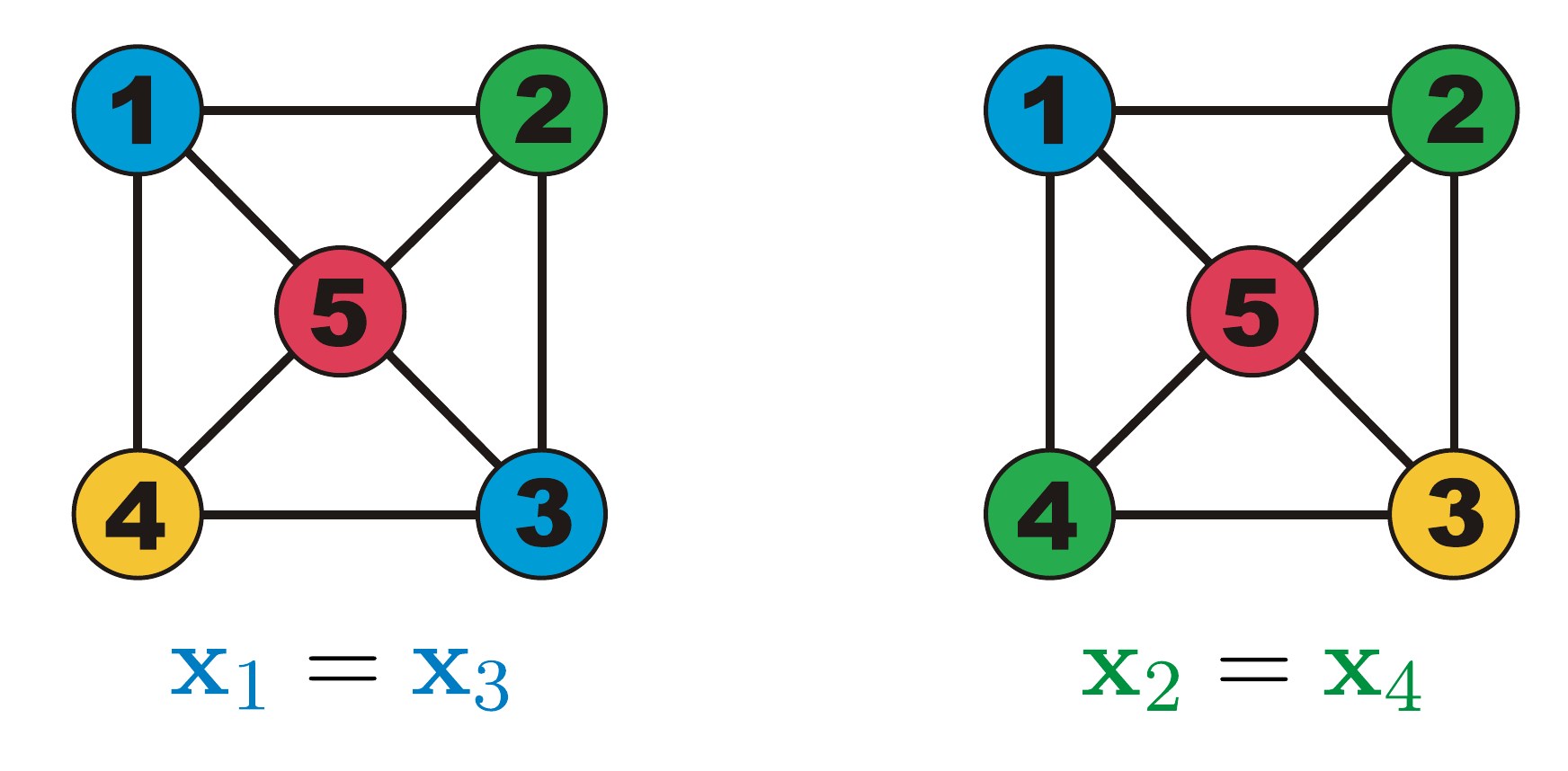}} & $\begin{pmatrix*}[r] 3 &\sqrt{2} &\sqrt{2} & \sqrt{2}& 0\\ \sqrt{2} & 3 &0 & -1 & 0  \\ \sqrt{2} &0 &3 & -1 & 0\\\sqrt{2} & -1 &-1 & 4 & 0  \\ 0 & 0 & 0 & 0 & 3\end{pmatrix*}$\\
\hline
A5 & $\mathbbm{1}$ & \raisebox{-.5\height}{\includegraphics[height=2.3cm]{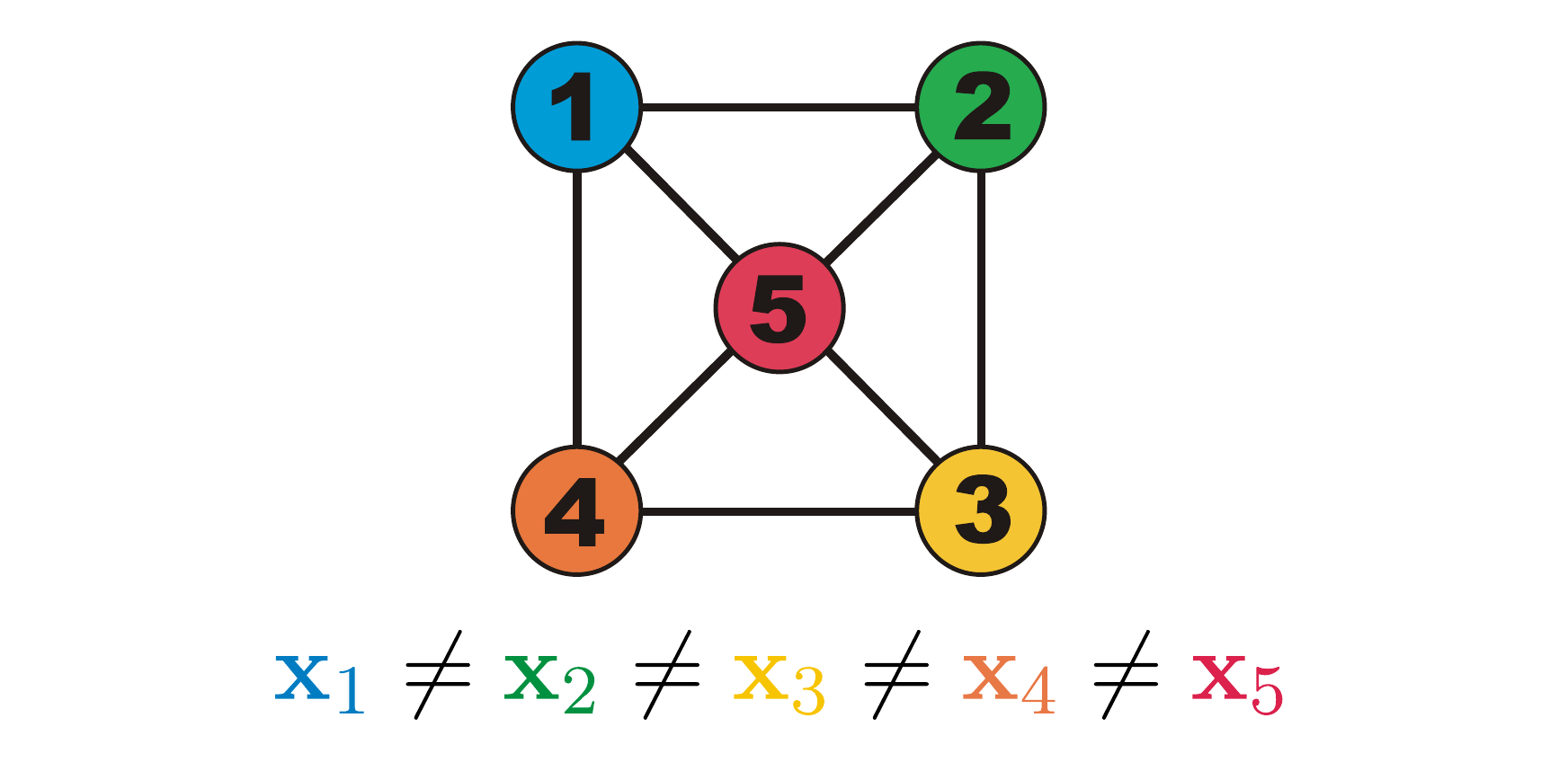}} & $\mathcal{G}$\\
\hline
\end{tabular}\caption{Cluster states due to symmetry in a network with connectivity given by \eqref{eq:Laplacian5}. Network symmetry group $\Gamma \cong \mathbf{D}_4$ gives clusters $\{1, 2, 3, 4\}$, $\{5\}$. Other potential cluster states are given by isotropy subgroups of $\Gamma$. Conjugate isotropy subgroups (such as the pairs in $A2$ and $A4$) give cluster states with identical existence and stability properties. Here we have used identical labelling for the cluster states to that given in \cite{Sorrentino2016}. We also give for each type of cluster state a corresponding block diagonalisation of the connectivity matrix, $\mathcal{G}^\prime= Q \mathcal{G}Q^{-1}$ for an appropriate transformation matrix $Q$.
\label{tab:clusters5}}
\end{center}
\end{table*}

%%%%%%%%%%%%%%%%%%%%%%%%%%%%%%%%%%%%%%%%%%%%%%%%%%%%%%%%%%%%%%%%%%%%%%%%%%%%%%%%%%%%%%%%%%%%%%%%%%%%%%%%%%%%%%%%%%
%%%%%%%%%%%%%%%%%%%%%%%%%%%%%%%%%%%%%%%%%%%%%%%%%%%%%%%%%%%%%%%%%%%%%%%%%%%%%%%%%%%%%%%%%%%%%%%%%%%%%%%%%%%%%%%%%%

Sorrentino \textit{et al}. \cite{Sorrentino2016} provide an algorithm for computing the subgroups of the network symmetry group $\Gamma$ which give rise to cluster states: Suppose that $S=\{ p_1, \ldots, p_{k}\}$ is a set of permutations which generates $\Gamma$. Partition $S$ into subsets $S= S_1 \cup \cdots \cup S_\nu$ such that the set of vertices moved by permutations in $S_i$ is disjoint from the set of vertices moved by permutations in $S_j$ for $i\neq j$, $i,j =1, \ldots \nu$, and each $S_i$ cannot itself be partitioned in this way. If we let $\mathcal{H}_i$ denote the subgroup generated by $S_i$ then we obtain the geometric decomposition of the group $\Gamma$ into the direct product of subgroups \[ \Gamma = \mathcal{H}_1 \times \cdots \times \mathcal{H}_\nu,\] \cite{MacArthur2008}. Note that for each node which is not permuted by any $p_i \in S$, we include a factor of $\mathcal{H}_i$ in the geometric decomposition where $S_i$ contains only the identity. Thus for the five-node example above $\Gamma = \mathcal{H}_1 \times \mathcal{H}_2$ where $S_1 = \{ \rho, \pi\}$ and $S_2 = \{ e\}$. One can use the geometric decomposition to determine all subgroups of $\Gamma$ which may give rise to cluster states by considering all groups of the type $\Sigma = \mathcal{K}_1 \times \cdots \times \mathcal{K}_\nu$ where $\mathcal{K}_i$ is a subgroup of $\mathcal{H}_i$. Efficient computational algebra routines \cite{GAP4, sagemath} are available which will compute the geometric decomposition and subgroups of a given arbitrary permutation group, so the computations can be carried out even for large networks with tens of thousands of nodes \cite{MacArthur2008}.

\subsection{\label{sec:ClustersLaplacian}Cluster states from Laplacian coupling}

In the general theory of patterns of synchrony in coupled cell systems \cite{Golubitsky2016}, cluster states correspond to \textit{balanced colourings} of the network nodes (where synchronous cells have equivalent inputs accounting for different types of cells and inputs and also self-coupling). Networks whose topology does not have any symmetries can support balanced colourings, while networks with symmetry can have additional balanced colourings alongside those which result from symmetry. When we consider networks with Laplacian coupling dynamics, as in \eqref{network1}, we expect to see additional potential cluster states resulting from balanced colourings which arise from the self-coupling. These additional states will include global synchronisation which is guaranteed to exist in networks of identical oscillators with Laplacian coupling.

An algorithm for computing these additional potential cluster states is given in \cite{Sorrentino2016} which uses as building blocks the clusters which can be found from network symmetries. First choose an isotropy subgroup $\Sigma = \mathcal{K}_1 \times \cdots \times \mathcal{K}_\nu \subseteq \Gamma$. This gives a partition of the nodes into clusters (a potential cluster state from symmetry). Next, merge some of these clusters as possibilities for new cluster states before finally checking if the merged cluster state (which cannot be one of those which arises due to symmetry) is dynamically valid. By this we mean that if we set all of the $\vec{x}_i$ to be equal for nodes in the merged clusters, the equations of motion are consistent. This checking for consistency in the network equations can be done by eye for small networks. However, a better way to complete this final step for large networks is to automate the checking process using computational algebra tools \cite{GAP4, sagemath}. To do this we use  the following steps which are described in greater detail in \cite{Sorrentino2016}. If the clusters we wish to merge are synchronised then the dynamics are equivalent to those which would be obtained if nodes in the same merged cluster were not connected (since the feedback term, $\vec{H}(\vec{x}_i)$,  will cancel coupling terms of nodes from the same merged cluster). This gives a dynamically equivalent graph Laplacian matrix $\overline{\mathcal{G}}$ which is the original graph Laplacian matrix $\mathcal{G}$ with off-diagonal entries between nodes in the same merged cluster set to $0$ and the diagonals set to the negative of the new row sums. If, when we compute the symmetry group $\overline{\Gamma}$ of the network with this graph Laplacian, one of its isotropy subgroups corresponds to our merged cluster state then the dynamics are flow-invariant and our merged cluster state is a dynamically valid synchronised state. Sorrentino \textit{et al}. \cite{Sorrentino2016} call the cluster states found in this way \textit{Laplacian} clusters. Note that there is no need to check mergings between subgroups of the same group $\mathcal{H}_i$.

\medskip
\noindent{\textbf{Example: Laplacian clusters}}

\smallskip
\noindent
For the same five-node network as discussed previously, we can see that if we take isotropy subgroup $\Sigma = \mathbb{Z}_2(\rho^2)$ which gives the clusters $\{ 1,3\}$, $\{ 2, 4\}$ and $\{5\}$ then one potential Laplacian merged cluster state may be $\{ 1,3,5\}$, $\{2,4\}$. This merging gives new graph Laplacian matrix
 \[ \overline{\mathcal{G}} = \begin{pmatrix*}[r] 2 &-1 &0 & -1 & 0\\ -1 & 3 &-1 & 0 & -1  \\ 0 &-1 &2 & -1 & 0\\-1 & 0 &-1 & 3 & -1  \\ 0 & -1 & 0 & -1 & 2\end{pmatrix*}.\] A network with this connectivity matrix has symmetry group $\overline{\Gamma}$ generated by $S= \{ (13), (35), (24)\}$. Thus the isotropy subgroup $\overline{\Gamma}$ gives the clusters $\{1,3,5\}$, $\{2,4\}$ and therefore this cluster state is dynamically valid. It is called $L3$ in Table \ref{tab:Lclusters5} (along with its conjugate cluster state).

Another potential Laplacian merged cluster state may be $\{ 1,2,5\}$, $\{3,4\}$ from merging clusters from the symmetry cluster state $A2$ (see Table \ref{tab:clusters5}). In this case the merging gives a new graph Laplacian matrix \[ \overline{\mathcal{G}} = \begin{pmatrix*}[r] 1 &0 &0 & -1 & 0\\ 0 & 1 &-1 & 0 & 0  \\ 0 &-1 &2 &0 & -1\\-1 & 0 &0 & 2 & -1  \\ 0 & 0 & -1 & -1 & 2\end{pmatrix*},\] which corresponds to a network with symmetry group $\overline{\Gamma}$ generated by $S= \{ (12)(34)\}$. Here the only isotropy subgroups are $\overline{\Gamma}$ which gives the cluster state $A2$ and the trivial subgroup which gives the cluster state $A5$. Thus this merging of clusters does not give a dynamically valid state.

Table \ref{tab:Lclusters5} gives a list of all dynamically viable Laplacian cluster states for the five-node network.  \hfill $\Diamond$
\medskip

 \begin{table*}[hbtp]
\begin{center}
\begin{tabular}{|c|c|c|}
\hline
Label  & Cluster state & $\mathcal{G}^{\prime\prime}$\\
\hline
L1 & \raisebox{-.5\height}{\includegraphics[height=2.3cm]{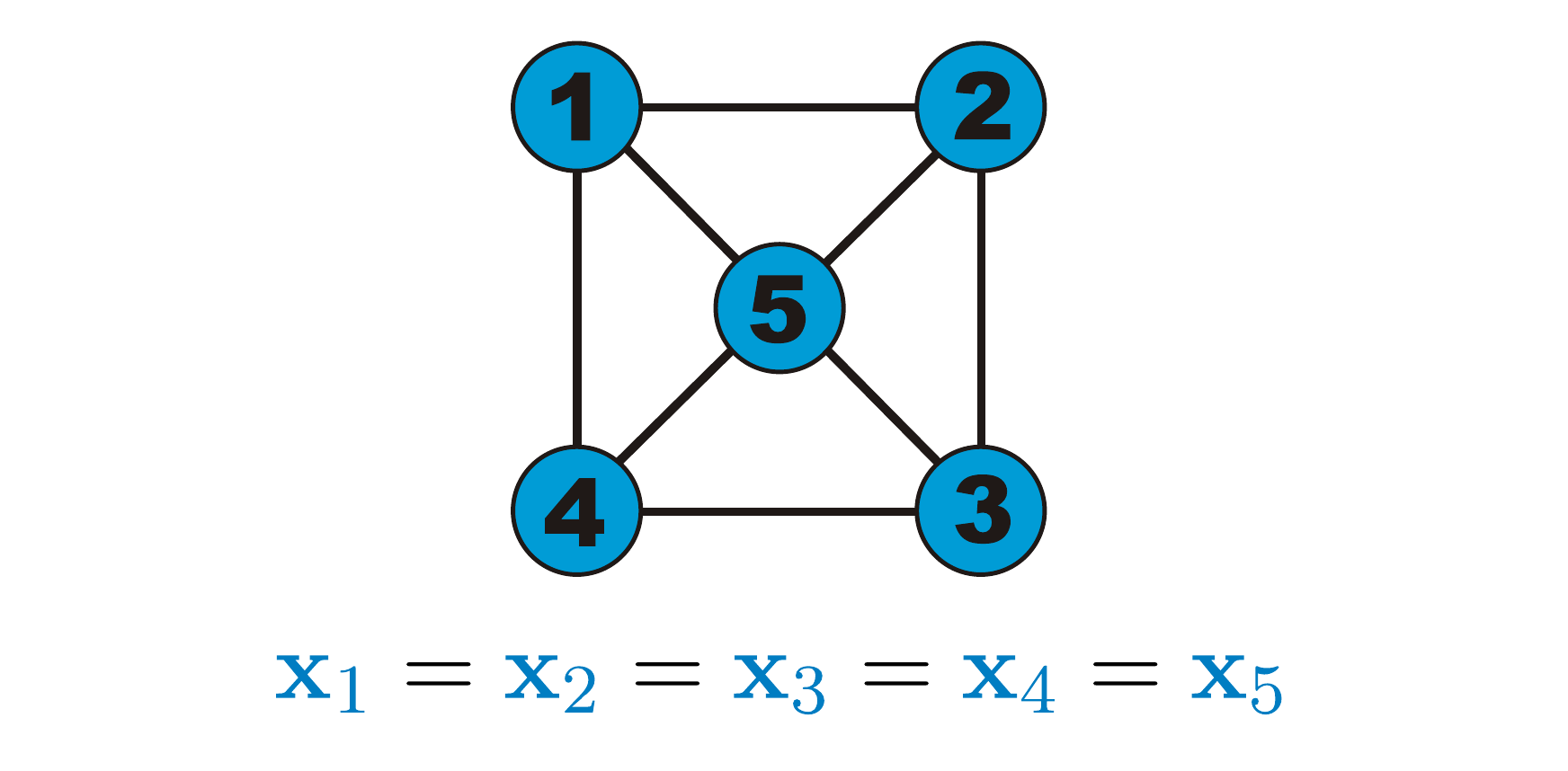}} & $\begin{pmatrix*}[r] 0 &0 &0 & 0 & 0\\ 0 & 3 &0 & 0 & 0  \\ 0 &0 &5 & 0 & 0\\0 & 0 &0 & 3 & 0  \\ 0 & 0 & 0 & 0 & 5\end{pmatrix*}$\\
\hline
L3 & \raisebox{-.5\height}{\includegraphics[height=2.3cm]{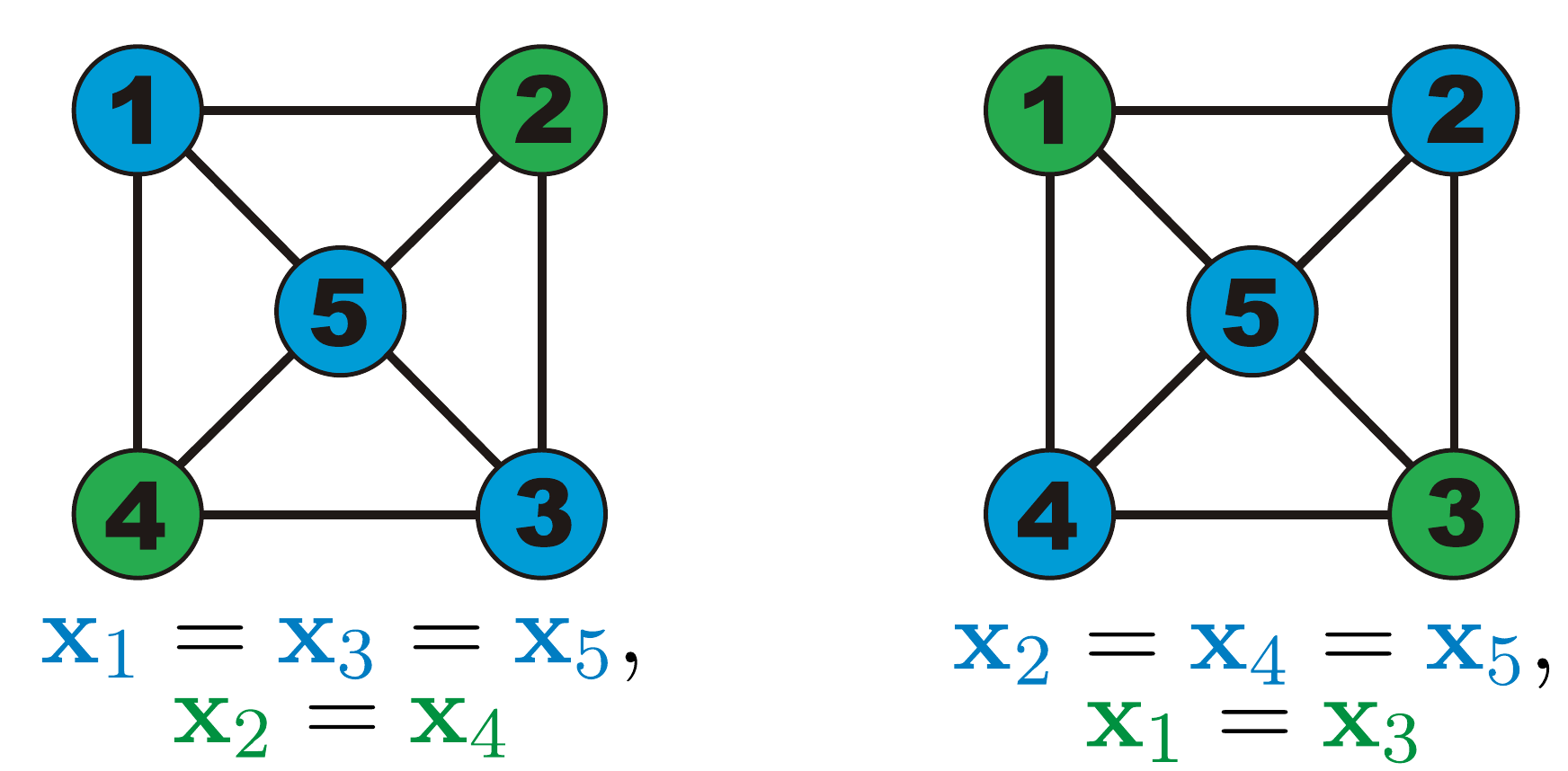}} & $\begin{pmatrix*}[r] 3 &-\sqrt{6} &0 & 0 & 0\\ -\sqrt{6} & 2 &0 & 0 & 0  \\ 0 &0 &5 & 0 & 0\\0 & 0 &0 & 3 & 0  \\ 0 & 0 & 0 & 0 & 3\end{pmatrix*}$\\
\hline
L4 & \raisebox{-.5\height}{\includegraphics[height=2.3cm]{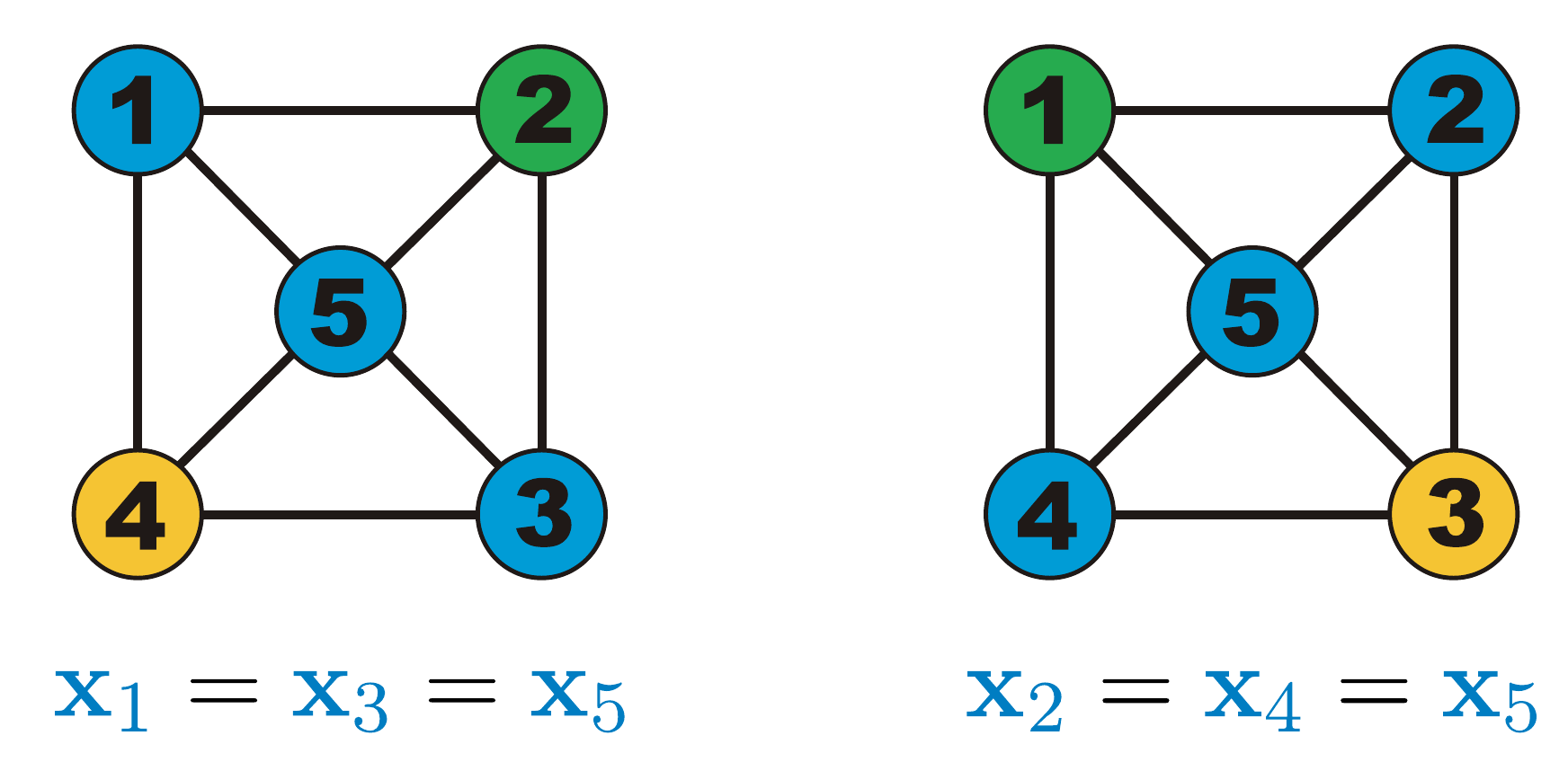}} & $\begin{pmatrix*}[r] 3 &-\sqrt{3} &0 & 0 & 0\\ -\sqrt{3} & 2 &-\sqrt{3} & 0 & 0  \\ 0 &-\sqrt{3} &3 & 0 & 0\\0 & 0 &0 & 5 & 0  \\ 0 & 0 & 0 & 0 & 3\end{pmatrix*}$\\
\hline
\end{tabular}\caption{Cluster states due to Laplacian coupling in a network with connectivity given by \eqref{eq:Laplacian5}. Here we have used identical labelling for the cluster states to that given in \cite{Sorrentino2016}. We also give for each type of cluster state a corresponding block diagonalisation of the connectivity matrix, $\mathcal{G}^{\prime\prime}= \chi \mathcal{G} \chi^{-1}$ for an appropriate transformation matrix $\chi$.
\label{tab:Lclusters5}}
\end{center}
\end{table*}

For Laplacian coupled networks with symmetry, Sorrentino \textit{et al}. \cite{Sorrentino2016} argue that their algorithm for computing cluster states which do not arise directly from symmetry is much more computationally efficient than algorithms not based on symmetries which compute all balanced colourings for any network topology, in particular that of Kamei and Cock \cite{Kamei2013}. They are careful to point out that while the checking stage (computing geometric decompositions and subgroups) is computationally efficient, the number of possible cluster mergings which must be checked grows combinatorially with the number $\mu$ of cluster states from symmetry with an upper bound of $B_\mu$ the $\mu$-th Bell number. Thus the algorithm of Sorrentino \textit{et al}. \cite{Sorrentino2016} is much faster than that of Kamei and Cock \cite{Kamei2013} when the number of cluster states from symmetry is small compared to the size of the network, that is $\mu \ll N$. Both Sorrentino \textit{et al}. \cite{Sorrentino2016} and Kamei and Cock \cite{Kamei2013} note that finding all dynamically valid cluster states for networks with Laplacian coupling may be substantially more difficult than for symmetric networks without self-coupling where connectivity is described by an adjacency matrix and all cluster states arise due to network symmetries.

\subsection{Stability of cluster states}

The methods outlined in sections \ref{sec:ClustersSymmetry} and \ref{sec:ClustersLaplacian} provide a catalogue of possible cluster states which may exist within a given network with graph Laplacian coupling \eqref{network1}. In applications, which of these states may exist and be stable will depend on the particular choices we make for the local dynamics $\vec{F}$ and the coupling function $\vec{H}$ as well as the global coupling strength $\sigma$. The presence of symmetry within the system imposes constraints on the form of the Jacobian matrix which can be used to greatly simplify stability calculations. Here we briefly review well-established methods for stability calculations within symmetric systems which apply to the cluster states which arise from network symmetries \cite{Golubitsky1988}. We also summarise the results of Sorrentino \textit{et al}. \cite{Sorrentino2016} which extend these techniques to Laplacian cluster states.

First consider a periodic cluster state arising from network symmetry that has isotropy $\Sigma \subseteq \Gamma$. The fixed-point subspace of this subgroup is $\Upsilon = \mbox{Fix}(\Sigma)$ which is the synchrony subspace for the cluster state. The state consists of $M$ clusters $\mathcal{C}_k$, $k=1, \ldots, M$, where $M = \operatorname{dim}( \operatorname{Fix} (\Sigma))= \dim(\Upsilon)$. Letting $\vec{s}_k(t)$ denote the synchronised state of nodes in cluster $\mathcal{C}_k$, and using the notation of \S \ref{sec:MSF}, the variational equation of \eqref{network1} about the cluster state is
\begin{align}
\label{eq:ClusterVariational}
\FD{}{t} \vec{U} =  &\left[ \sum_{k=1}^M E^{(k)} \otimes D \vec{F}(\vec{s}_k(t)) \right. \nonumber \\
& \left. - \sigma \sum_{k=1}^M \left( \mathcal{G} E^{(k)}\right) \otimes D\vec{H}(\vec{s}_k(t)) \right ] \vec{U} ,
\end{align}
where $E^{(k)}$ is the diagonal $N \times N$ matrix such that $E_{ii}^{(k)} = 1$ if $i \in \mathcal{C}_k$ and $E_{ii}^{(k)} = 0$ otherwise. To determine the stability of the periodic cluster state we need to compute the Floquet exponents of \eqref{eq:ClusterVariational}. This task can be greatly simplified due to the fact that the system of variational equations can be block-diagonalised using the symmetries present in the system.

The node space can be decomposed into a number of irreducible representations of the isotropy subgroup $\Sigma$. Some of these subspaces will be isomorphic to each other and we combine these to obtain the isotypic components of the node space. Each isotypic component is invariant under the variational equation \eqref{eq:ClusterVariational} so Floquet exponents can be found by considering the restriction of the equations to each isotypic component. Thus the decomposition puts the variational equations into block-diagonal form and we then compute Floquet exponents for each block to determine stability of the cluster state. The process of isotypic decomposition and its use in stability computations is derived in detail in \cite{Golubitsky1988}. Pecora \textit{et al}. \cite{Pecora2014} give an explicit algorithm to determine the isotypic decomposition for a given cluster state from symmetry and compute the transformation matrix $Q$ such that $\mathcal{G}^\prime = Q \mathcal{G} Q^{-1}$ is block diagonal. Applying this transformation to the variational equation \eqref{eq:ClusterVariational} we obtain block diagonal system of equations
\begin{align}
\label{eq:blockVariational}
\FD{}{t}{\vec{V}} = & \left[\sum_{k=1}^M J^{(k)} \otimes D \vec{F}(\vec{s}_k(t)) \right. \nonumber \\
& \left.- \sigma \sum_{k=1}^M \left( \mathcal{G}^\prime J^{(k)}\right) \otimes D\vec{H}(\vec{s}_k(t)) \right] \vec{V} ,
\end{align}
where $\vec{V}(t) = (Q \otimes I_m) \vec{U}(t)$ and $J^{(k)} = Q E^{(k)} Q^{-1}$. In Table \ref{tab:clusters5} we give the block-diagonalised graph Laplacian matrix for each of the cluster states from symmetry. An important point to note is that the isotypic component of the trivial representation is $\mbox{Fix}(\Sigma) = \Upsilon$, the synchronisation manifold. The block corresponding to this component is the $M \times M$ block which appears in the top left in each of the examples in Table \ref{tab:clusters5}. This block corresponds to perturbations within the synchronisation manifold and will always have a Floquet exponent equal to $0$.  The remaining blocks correspond to the isotypic components of other irreducible representations of $\Sigma$. When the node space representation contains $l\geq 1$ isomorphic copies of a particular irreducible representation this will result in a block of dimension $l$. For example, see the $2 \times 2$ block in the block-diagonalisation for cluster state $A2$ in Table \ref{tab:clusters5}. These blocks represent perturbations transverse to the synchronisation manifold and Floquet multipliers from these blocks will determine stability under synchrony-breaking perturbations. For the cluster state to be stable all Floquet exponents (except the one which is always $0$) must have negative real part.

In the case of a periodic Laplacian cluster state, the synchronisation manifold is an invariant subspace, but it is not the fixed-point subspace of any subgroup of $\Gamma$. However, we can still find a block-diagonalisation of $\mathcal{G}$ which has in the top left a block which corresponds to perturbations within the synchronisation manifold. We do this following the algorithm of Sorrentino \textit{et al}. \cite{Sorrentino2016}. Suppose that we start with a cluster state from symmetry with isotropy $\Sigma$ which has $M$ clusters and whose variational equations are block diagonalised by the matrix $Q$. Now suppose that a Laplacian cluster state is obtained by merging together two of the clusters in this state. The dimension of the synchronisation manifold decreases by one, while the dimension of the transverse manifold increases by one. New coordinates on the synchronisation manifold are obtained by transforming the new synchronisation vector in the node space coordinates (this has a 1 in the position of every node in the new merged cluster and $0$s elsewhere) into the coordinates of the block-diagonalisation of the cluster state with isotropy $\Sigma$. The orthogonal complement provides the new transverse direction. Normalising the resulting vectors and entering them as rows of an orthogonal matrix $Q^\prime$ whose other rows have $Q^\prime_{ij}=\delta_{ij}$, we find that the matrix $\chi=Q^\prime Q$ block-diagonalises $\mathcal{G}$ to a matrix $\mathcal{G}^{\prime \prime}$ which has top left block of dimension $(M-1) \times (M-1)$. Thus the transformation matrix $\chi$  will block diagonalise the variational equations for the Laplacian cluster state so we may more easily determine the $m(M-1)$ Floquet exponents within the synchronisation manifold and the $m(M+1)$ transverse Floquet exponents.

\medskip
\noindent{\textbf{Example: Block diagonalisation}}

\smallskip
\noindent
We revisit the example five-node network of Sorrentino \textit{et al}. \cite{Sorrentino2016}. The block-diagonalisations for the cluster state $A3$ and the merging which gives cluster state $L3$ are covered in \cite{Sorrentino2016}. Here we consider the cluster states $A4$ and $L4$ for illustration of the algorithms outlined above. The cluster state $A4$ has (choosing one option from the conjugacy class as they have identical existence and stability properties) isotropy subgroup $\mathbb{Z}_2(\pi)$. This group has only two irreducible representations, the trivial representation and the representation where $\pi$ acts as multiplication by $-1$. The node space representation contains 4 isomorphic copies of the trivial representation and one copy of the other representation. Following the algorithm given explicitly in \cite{Pecora2014} we obtain
\[ Q = \frac{1}{2} \begin{pmatrix*}[r] 2 &0&0 & 0& 0\\ 0 & -\sqrt{2} &0 & -\sqrt{2} & 0  \\ 0 &0 &2 & 0 & 0\\0& 0 &0 & 0 & 2  \\ 0 & -\sqrt{2}  & 0 & \sqrt{2}  & 0\end{pmatrix*}, \] and $\mathcal{G}^\prime$ as in Table \ref{tab:clusters5}.

Now consider the Laplacian cluster state $L4$ which is found by merging $\{2,4\}$ with $\{5\}$ to give a two cluster state $\{2,4,5\}$, $\{1,3\}$. The new synchronisation vector in node coordinates will be $\vec{v}_1 = (0,1,0,1,1)^T$ and new transverse vector is $\vec{v_2} = (0,1,0,1, -2)^T$. Normalising $Q\vec{v}_1$ and $Q\vec{v}_2$ we see that we must put these vectors as rows 2 and 4 of the orthogonal matrix $Q^\prime$ to obtain
\[ Q^\prime = \frac{1}{3} \begin{pmatrix*}[r] 3 &0&0 & 0& 0\\ 0 & -\sqrt{6} &0 & \sqrt{3} & 0  \\ 0 &0 &3 & 0 & 0\\0& -\sqrt{3} &0 & -\sqrt{6} & 0  \\ 0 & 0  & 0 &0  & 3\end{pmatrix*},\] and $\mathcal{G}^{\prime\prime}= Q^\prime \mathcal{G}^\prime (Q^\prime)^{-1}$ as in Table \ref{tab:Lclusters5}. \hfill $\Diamond$

\section{\label{sec:Nonsmooth}Orbits and variational equations for nonsmooth systems}

Using the algorithms in section \ref{sec:Cluster} one can, given a network structure, compute the catalogue of cluster states (from symmetry and Laplacian coupling) which may be observed for given local dynamics $\vec{F}$ and interaction function $\vec{H}$. It is also possible, to block-diagonalise the variational equations needed to compute the stability of a given cluster state. However in practice, to determine which of the cluster states are stable and to carry out any kind of bifurcation analysis, we observe from \eqref{eq:blockVariational} that we will require closed form solutions for the periodic cluster state. For many of the situations in which one would wish to apply the framework, this information simply is not available. Therefore, we turn our attention to networks of PWL oscillators where, as we will show here, it is relatively straightforward to construct the periodic orbits for the cluster state and additionally apply the required modifications to the Floquet theory to account for the lack of smoothness of the dynamics.

Here we demonstrate how the computations of periodic cluster states and their stability may be carried out for a $2$-cluster state in a network with $2$-dimensional local dynamics and one switching plane. The method outlined here can easily be extended to larger numbers of clusters and more complex local PWL dynamics (e.g. with more switching planes and jump discontinuities) and we consider such examples in section \ref{sec:Applications}.

\begin{figure}[t]
\begin{center}
\includegraphics[width=6cm]{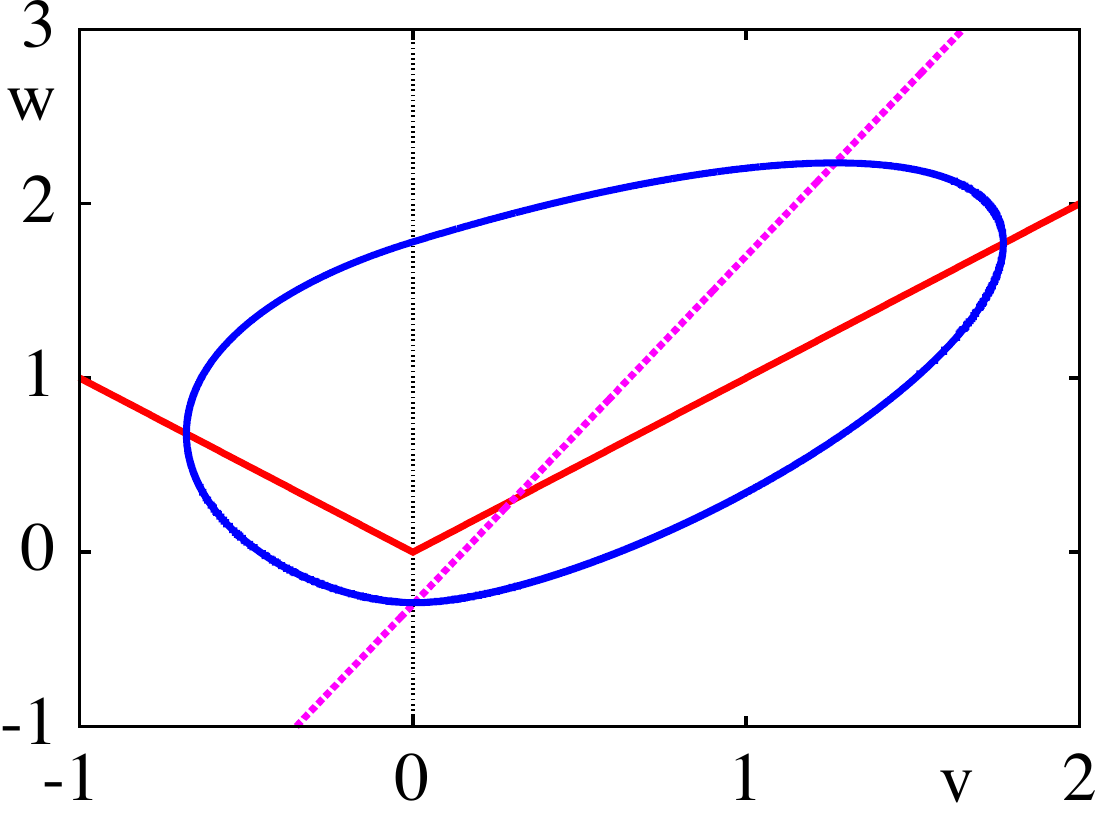}
\caption{Phase plane of the absolute model with $\bar{v} = 0.1$, $\bar{w}=-0.1$ and $g=0.5$. The $v$-nullcine is in red (solid) and the $w$-nullcline is in green (dashed). The unstable equilibrium point lies in $S_R$ and a stable periodic orbit is shown in blue.}
\label{Fig:Absolute}
\end{center}
\end{figure}

Consider the specific case of a two cluster state $L3$ in the five-node example network with graph Laplacian given by \eqref{eq:Laplacian5}. Let each of the nodes have two-dimensional local dynamics given by the \textit{absolute} model, which is named because it has one nullcline described by the absolute value function (and see \cite{Coombes2011a} for a further discussion), where
\begin{equation}
A_L=\begin{bmatrix}
-1 & -1\\
1 & -g
\end{bmatrix}, \ \
A_R=\begin{bmatrix}
1& -1\\
1 & -g
\end{bmatrix}, \ \
\vec{c}_L=\begin{bmatrix}
 0\\
 g\bar{w}-\bar{v}
\end{bmatrix}= \vec{c}_R,
\nonumber
\end{equation} with $0<g<1$, $g\bar{w}<\bar{v}$ and indicator function $h(\vec{z}; 0) = v$. This model supports a non-smooth Andronov--Hopf bifurcation which occurs when the equilibrium moves from $v<0$ to $v>0$ crossing the switching manifold. When the (unstable) equilibrium lies in $S_R$, a periodic orbit for the single node dynamics can be constructed by connecting two trajectories, one in each of the half planes. Starting from initial data $\vec{z}(0) = (0, w(0))^\mathsf{T}$, this requires that we simultaneously solve the three nonlinear equations $v(T_R)=0$, $v(\Delta)=0$, $w(\Delta)=w(0)$ where $\Delta$ is the period of the orbit, $T_R$ is the time of flight in $S_R$ and thus $\Delta>T_R>0$. We do this numerically using Matlab's fsolve routine.  The periodic orbit can be shown to be stable with Floquet exponent \cite{Coombes2016}
\[ r = -g - \frac{\Delta-2T_R}{\Delta}.\] The nullclines for the model together with the stable periodic orbit are shown in Fig.~\ref{Fig:Absolute}.

%%%%%%%%%%%%%%%%%%%%%%%%%%%%%%%%%%%%%%%%%%%%%%%%%%%%%%%%%%%%%%%%%%%%%%%%%%%%%%%%%%%%%%%%%%%%%%%%%%%%%%%%%%%%%%%%%%
%%%%%%%%%%%%%%%%%%%%%%%%%%%%%%%%%%%%%%%%%%%%%%%%%%%%%%%%%%%%%%%%%%%%%%%%%%%%%%%%%%%%%%%%%%%%%%%%%%%%%%%%%%%%%%%%%%

\begin{figure}[t]
\begin{center}
\includegraphics[width=6cm]{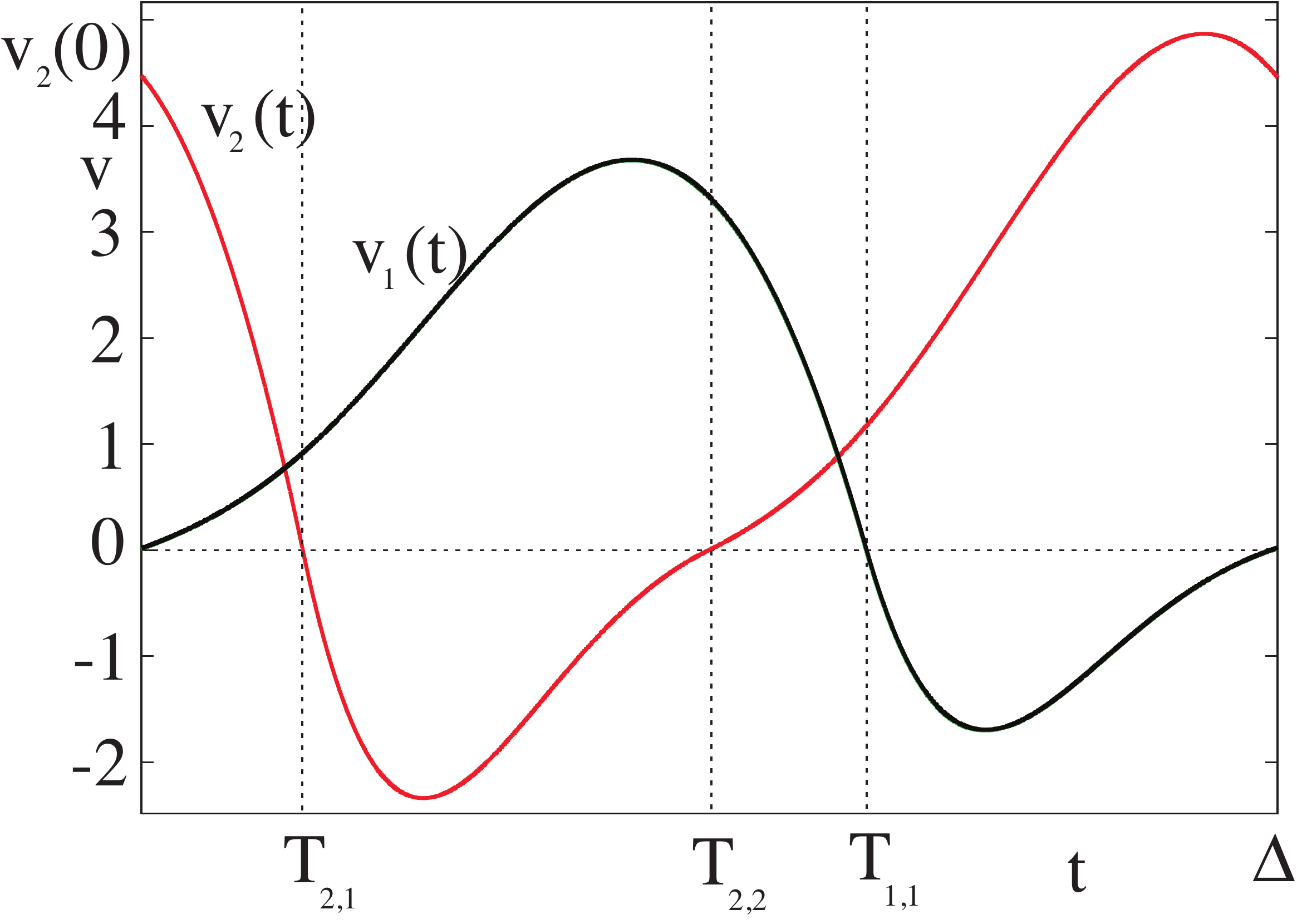}
\caption{The $v$ components of the orbits $\vec{s}_1$ and $\vec{s}_2$ over one period. Seven nonlinear algebraic equations must be solved to find the unknown initial data $w_1(0), v_2(0), w_2(0)$ and switching times $T_{1,1}$, $T_{2,1}$, $T_{2,2}$ and $T_{1,2}= \Delta$ as shown.}
\label{Fig:switching}
\end{center}
\end{figure}

Now couple these nodes together into a network described by \eqref{network1} where $\mathcal{G}$ is given by \eqref{eq:Laplacian5} and the coupling function is given by $\vec{H}(\vec{z}) = (v,0)^\mathsf{T}$ (diffusive coupling). Numerical simulations of the network equations with the parameter values of Fig.~\ref{Fig:Absolute} indicate that a periodic $L3$ cluster state exists with coupling strength $\sigma=-0.03$. On this cluster state $\vec{z}_1 = \vec{z}_3 =\vec{z}_5 = \vec{s}_1$ and $\vec{z}_2 = \vec{z}_4= \vec{s}_2$. (Recall that the conjugate cluster state where $\vec{z}_1 = \vec{z}_3$ and $\vec{z}_2 = \vec{z}_4= \vec{z}_5$ has the same existence and stability properties.) On this invariant subspace the equations have the form $\dot{\vec{s}} = A_{\mu_1, \mu_2} \vec{s} + \vec{c}_{\mu_1, \mu_2}$ where $\vec{s} = (\vec{s}_1, \vec{s}_2)^\mathsf{T}$,
\[ A_{\mu_1, \mu_2} = \begin{bmatrix}
A_{\mu_1} - 2 \sigma D\vec{H} & 2 \sigma D\vec{H}\\
3 \sigma D\vec{H} & A_{\mu_2} -3 \sigma D\vec{H}
\end{bmatrix}, \ \
\vec{c}_{\mu_1, \mu_2} = \begin{bmatrix}
\vec{c}_{\mu_1} \\
\vec{c}_{\mu_2}
\end{bmatrix},
\] and \[\mu_i=  \begin{cases} L, & v_i<0\\ R, & v_i>0\end{cases} .\]
This is simply a $4$-dimensional PWL system with two switching planes, $v_1=0$ and $v_2=0$. The periodic orbit on the $4$-dimensional synchronous manifold can again be constructed by connecting together trajectories which satisfy the linear ODEs in each region across the switching planes. Starting from initial data $\vec{s}(0) = (0, w_1(0), v_2(0), w_2(0))^\mathsf{T}$, we now have to solve a system of seven nonlinear algebraic equations for $w_1(0), v_2(0), w_2(0)$ and the four switching times $T_{1,1}$, $T_{2,1}$, $T_{2,2}$ and $T_{1,2}= \Delta$. See Fig.~\ref{Fig:switching}.

%%%%%%%%%%%%%%%%%%%%%%%%%%%%%%%%%%%%%%%%%%%%%%%%%%%%%%%%%%%%%%%%%%%%%%%%%%%%%%%%%%%%%%%%%%%%%%%%%%%%%%%%%%%%%%%%%%
%%%%%%%%%%%%%%%%%%%%%%%%%%%%%%%%%%%%%%%%%%%%%%%%%%%%%%%%%%%%%%%%%%%%%%%%%%%%%%%%%%%%%%%%%%%%%%%%%%%%%%%%%%%%%%%%%%

The initial data and switching times can now be used to compute explicitly the Floquet multipliers for the periodic orbit, making use of the block-diagonalisation of the variational equations \eqref{eq:blockVariational}. Multipliers corresponding to perturbations within the synchronisation manifold can be computed without using the block-diagonalisation. We have
\[ \FD{}{t}\delta {\vec{s}} = A_{\mu_1, \mu_2}\delta\vec{s} , \] which can be solved using matrix exponentials, being careful to use saltation matrices to evolve perturbations through switching manifolds. After one period \[\delta\vec{s}(\Delta) = \Psi_s\delta\vec{s}(0) , \] where from Fig.~\ref{Fig:switching},
\begin{align}
\Psi_s = &K_{12} \mathcal{E}_{L,R}(\Delta- T_{1,1}) K_{11}  \mathcal{E}_{R,R}(T_{1,1}-T_{2,2}) \nonumber \\
& \times K_{22}  \mathcal{E}_{R,L}(T_{2,2}-T_{2,1})K_{21}  \mathcal{E}_{R,R}(T_{2,1}),
\end{align}
with saltation matrices
\[ K_{ij} = P_i \otimes K_i(T_{i,j}), \quad
P_1 = \begin{bmatrix}
 1& 0 \\
0 & 0
\end{bmatrix}, \quad P_2 = \begin{bmatrix}
0& 0 \\
0 & 1
\end{bmatrix},\]
\begin{equation} \label{eq:saltation} K_i(t) = \begin{bmatrix}
\dot{v}_i(t^+)/\dot{v}_i(t^-) & 0 \\
(\dot{w}_i(t^+)-\dot{w}_i(t^-))/\dot{v}_i(t^-) & 1
\end{bmatrix}, \quad i=1,2,\end{equation}
and \[ \mathcal{E}_{\mu_1, \mu_2}(t) = \e^{A_{\mu_1, \mu_2}t}, \qquad \mu_i \in \{ L, R\}.\]
The Floquet multipliers are the eigenvalues of $\Psi_s$. One of these eigenvalues will always be $1$ corresponding to perturbations along the periodic orbit.

In the directions transverse to the synchronisation manifold, the block-diagonalisation of $\mathcal{G}$ results in the following three decoupled Floquet problems:
\begin{align}
\dot{\vec{V}}_3 &= (D \vec{F}(\vec{s}_1) -3\sigma D \vec{H})\vec{V}_3,  \nonumber \\
\dot{\vec{V}}_4 &= (D \vec{F}(\vec{s}_2) -3\sigma D \vec{H})\vec{V}_4,  \nonumber \\
\dot{\vec{V}}_5 &= (D \vec{F}(\vec{s}_1) -5\sigma D \vec{H})\vec{V}_5,  \nonumber
\end{align} which again can be solved using matrix exponentials and saltation matrices: $\vec{V
}_i(\Delta) = \Psi_{t_i} \vec{V}_i(0)$ where
\begin{align}
\Psi_{t_3} &= K_1(\Delta) \mathcal{E}_L^3(\Delta-T_{1,1}) K_1(T_{1,1}) \mathcal{E}_R^{3}(T_{1,1})  ,  \nonumber \\
\Psi_{t_4} &=\mathcal{E}_R^3(\Delta- T_{2,2}) K_2(T_{2,2}) \mathcal{E}_L^3(T_{2,2}-T_{2,1})\nonumber \\ & \times K_2(T_{2,1}) \mathcal{E}_R^3(T_{2,1}) ,  \nonumber \\
\Psi_{t_5} &= K_1(\Delta) \mathcal{E}_L^5(\Delta-T_{1,1}) K_1(T_{1,1}) \mathcal{E}_R^5(T_{1,1}) , \nonumber
\end{align}
and \[\mathcal{E}_\mu^{\beta}(t) = \e^{(A_\mu-\beta\sigma D\vec{H})t}.\]
The eigenvalues of the $\Psi_{t_j}$, $j \in {3,4,5}$ give the Floquet multipliers for directions transverse to the synchronisation manifold.

Note that the change of basis from $\vec{U}$ coordinates to $\vec{V}$ coordinates has no effect on the action of the saltation matrices:
Recall $\vec{V} = (Q \otimes I_m) \vec{U}$. To progress $\vec{U}$ through a switch or discontinuity then $\vec{U}^+ = K \vec{U}^-$, where
\[ K = \sum_{k=1}^M E^{(k)} \otimes K_k.\]
Thus $\vec{V}^+ =\widehat{K}  \vec{V}$ where
\begin{align*}
\widehat{K} &= (T \otimes I_m) K (T\otimes I_m)^{-1}  = \sum_{k=1}^M (T E^{(k)} T^{-1} \otimes K_k) \\
& = \sum_{k=1}^M (J^{(k)} \otimes K_k).
\end{align*}

\begin{figure}[t]
\begin{center}
\includegraphics[width=8cm]{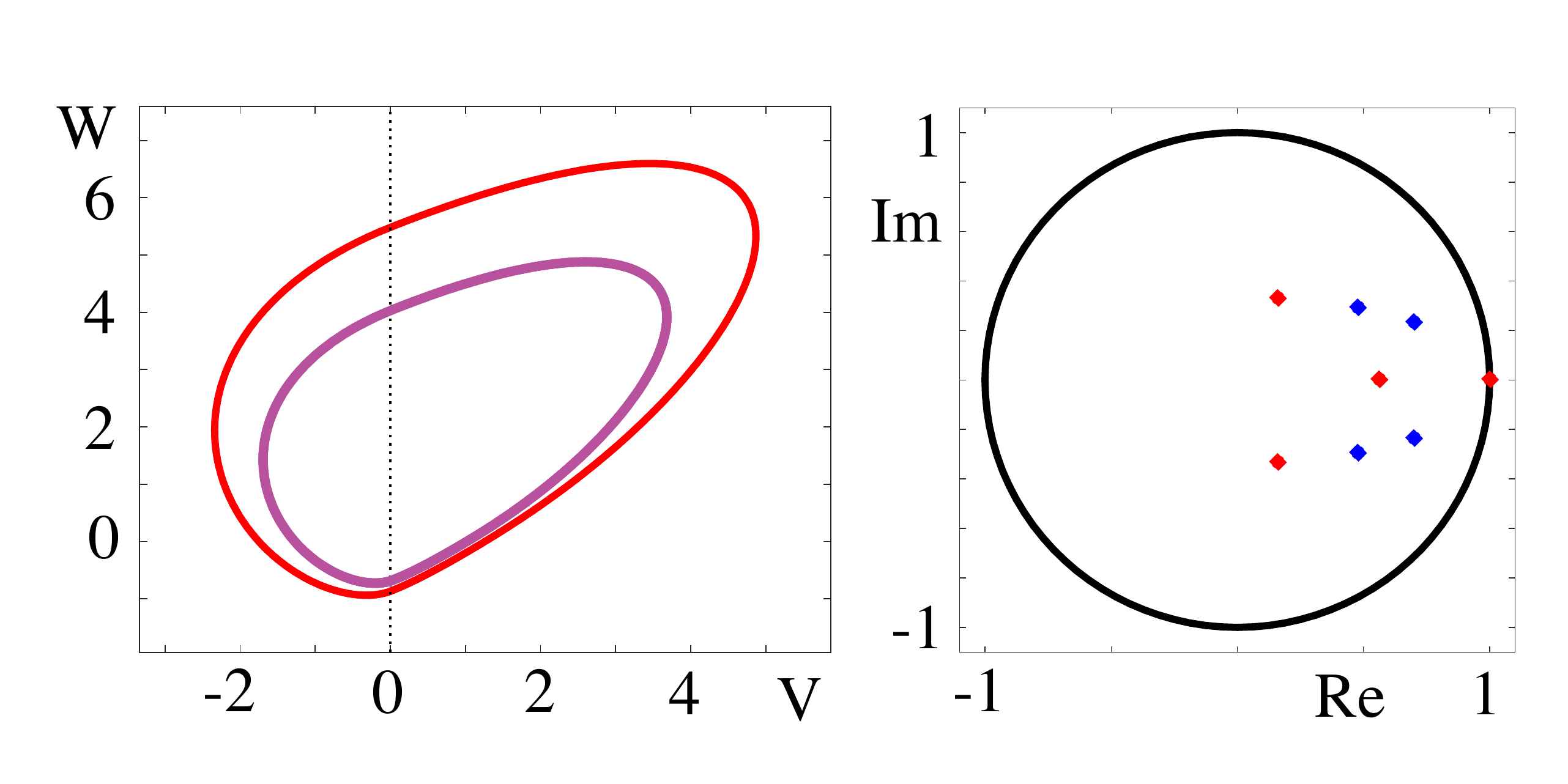}
\caption{Left:  The periodic orbits $\vec{s}_1$ (magenta) and $\vec{s}_2$ (red) for absolute model local dynamics with values of parameters as in Fig.~\ref{Fig:Absolute} and coupling strength $\sigma=-0.03$. The period $\Delta \approx 9.16$.
Right:  The Floquet multipliers for the periodic orbit lie inside the unit disc indicating that the orbit is stable. Floquet multipliers within the synchronisation manifold are shown in red, transverse directions are in blue.}
\label{Fig:stability}
\end{center}
\end{figure}

Since the vector field for the absolute model is continuous, in this case all saltation matrices are the identity. Figure \ref{Fig:stability} shows the Floquet multipliers for the choice of parameter values as in Fig.~\ref{Fig:Absolute}. As they all lie inside the unit disc (with the exception of the one which is forced to be $1$) the periodic orbit, as shown in the left hand panel of Fig.~\ref{Fig:stability} is stable.

%%%%%%%%%%%%%%%%%%%%%%%%%%%%%%%%%%%%%%%%%%%%%%%%%%%%%%%%%%%%%%%%%%%%%%%%%%%%%%%%%%%%%%%%%%%%%%%%%%%%%%%%%%%%%%%%%%
%%%%%%%%%%%%%%%%%%%%%%%%%%%%%%%%%%%%%%%%%%%%%%%%%%%%%%%%%%%%%%%%%%%%%%%%%%%%%%%%%%%%%%%%%%%%%%%%%%%%%%%%%%%%%%%%%%

\begin{figure*}[t]
\begin{center}
\includegraphics[width=16cm]{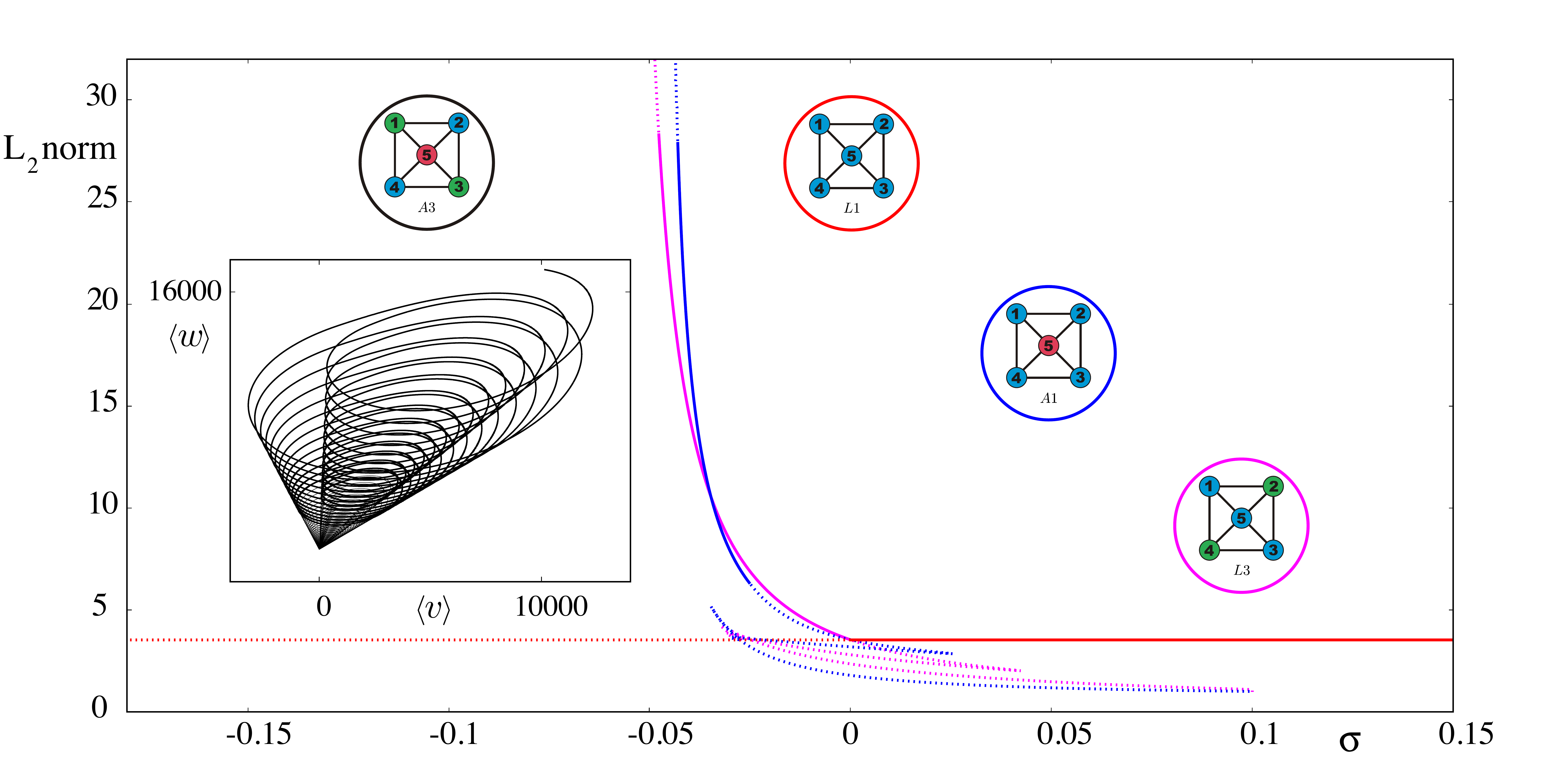}
\caption{Bifurcations between cluster states upon variation of coupling strength $\sigma$ in the five-node network of nodes with absolute model local dynamics (parameters as in Fig.~\ref{Fig:Absolute}). Stable periodic solutions are denoted by solid lines while unstable solutions are denoted by dotted lines. Each branch is coloured according to the circle containing the relevant cluster state. In the case of the $L3$ branch where a conjugate solution is expected with identical stability properties only one branch is shown. The inset shows the mean field dynamics, $(\langle v \rangle, \langle w \rangle ) = \sum_{i=1}^5 (v_i, w_i)/5$ , for $\sigma= -0.05$. The dynamics blow up in finite time whilst remaining as an $A3$ cluster state. This behaviour dominates for $\sigma \lesssim -0.0477$ where the $L3$ branch loses stability. Here all bifurcations from stable states are of tangent type where a Floquet multiplier passes through $+1$.}
\label{Fig:AbsoluteBifDiag}
\end{center}
\end{figure*}

We can locate bifurcations of the periodic orbit by determining when Floquet multipliers leave the unit disc (noting that the ordering of times which the trajectories cross the switching planes may also change as parameters are varied). Treating all cluster types we can build up bifurcation diagrams. For the absolute model with the parameters for the local dynamics as given in Fig.~\ref{Fig:Absolute}, bifurcations which can be found upon varying coupling strength $\sigma$ are shown in Fig.~\ref{Fig:AbsoluteBifDiag}. All bifurcations from stable states are of tangent type where a Floquet multiplier passes through $+1$.

%%%%%%%%%%%%%%%%%%%%%%%%%%%%%%%%%%%%%%%%%%%%%%%%%%%%%%%%%%%%%%%%%%%%%%%%%%%%%%%%%%%%%%%%%%%%%%%%%%%%%%%%%%%%%%%%%%
%%%%%%%%%%%%%%%%%%%%%%%%%%%%%%%%%%%%%%%%%%%%%%%%%%%%%%%%%%%%%%%%%%%%%%%%%%%%%%%%%%%%%%%%%%%%%%%%%%%%%%%%%%%%%%%%%%

\section{\label{sec:Applications}Further Applications and examples}

Extending the techniques illustrated in section \ref{sec:Nonsmooth}, here we discuss the cluster patterns and bifurcations which occur for different PWL local dynamics on the five node network \eqref{eq:Laplacian5}. We consider a diffusively coupled network of PWL Morris--Lecar oscillators (which has more than one switching plane) and the integrate-and-fire network model with event-driven synaptic coupling as discussed for the case of global synchrony in section \ref{sec:MSF-PWL} (which has jump discontinuities within the vector field). We show that each model of local dynamics gives a very different variety of cluster state dynamics demonstrating that emergent network behaviour depends as much on node dynamics as network structure.

\subsection{\label{sec:pwlml} Piecewise linear Morris--Lecar network}

\begin{figure}[t]
\begin{center}
\includegraphics[width=6cm]{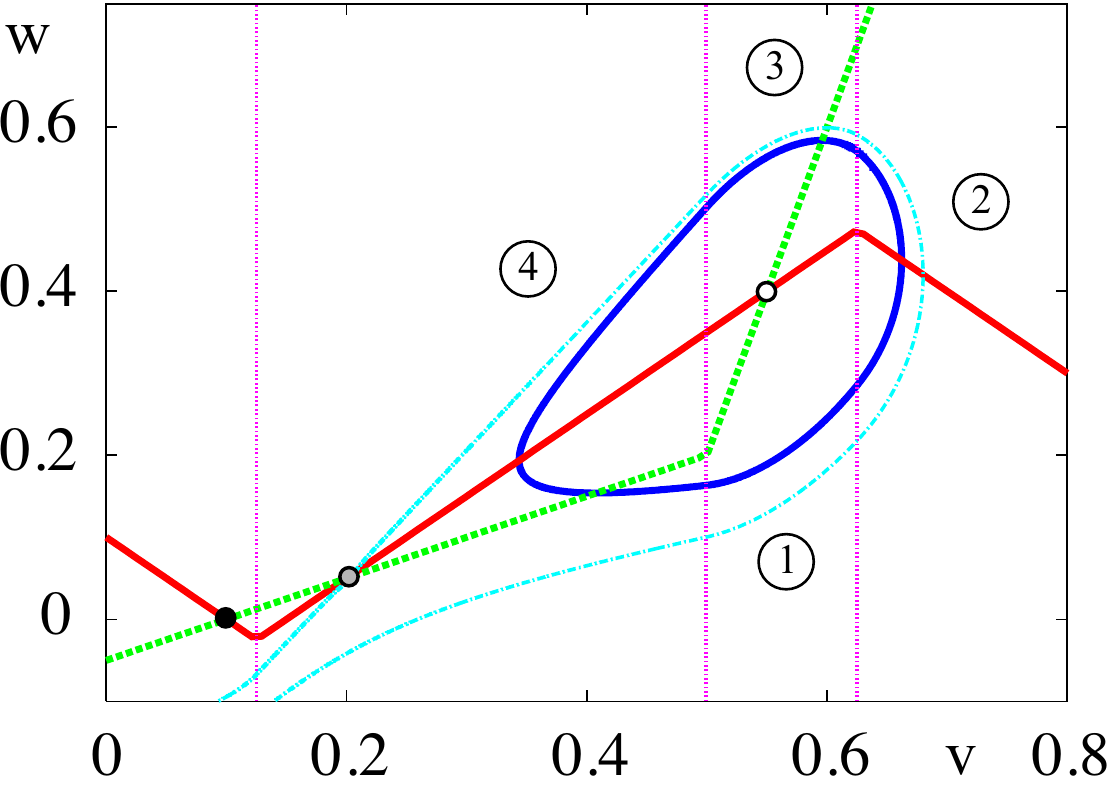}
\caption{The phase plane for the piecewise linear Morris--Lecar (PML) model with $\gamma_1$ = 2, $\gamma_2 = 0.25$, $C = 0.825$, $I = 0.1$, $a = 0.25$, $b = 0.5$ and $b^* = 0.2$. The $v$-nullcline is shown in red and the $w$-nullcline in magenta. The filled black circle indicates a stable fixed point, the grey filled circle a saddle and the filled white circle an unstable fixed point. The periodic orbit is shown in blue.
The pale blue line passing through the saddle (gray filled circle) is the separatrix between the stable fixed point (black filled circle) and the stable limit cycle (in blue).}
\label{Fig:PWLML}
\end{center}
\end{figure}

We now investigate dynamics of a network of planar PWL nodes with more than one switching plane. Consider the five-node network given by graph Laplacian coupling matrix \eqref{eq:Laplacian5} where the local dynamics $\vec{F}$ are given by a PWL reduction of the Morris--Lecar model (PML) \cite{Coombes2008}, and the coupling function is again given by $\vec{H}(\vec{z}) = (v,0)^\mathsf{T}$ (diffusive coupling). Figure \ref{Fig:PWLML} shows the phase plane for this model in the regime where there are type I oscillations (emerging with arbitrarily low frequency from a homoclinic bifurcation) and bistability between a periodic orbit and a stable fixed point.  The dynamics of the PWL model are defined by
\begin{align*}
C \dot{v} =  f(v)-w +I, \qquad
\dot{w} =  g(v,w),
\end{align*}
where
\begin{equation}
f(v) = \begin{cases}
-v, & v<a/2\\
v-a, & a/2 \leq v \leq (1+a)/2\\
1-v, & v > (1+a)/2
\end{cases} ,
\nonumber
\end{equation}
and
\begin{equation}
g(v,w) = \begin{cases}
( v-\gamma_1 w + b^*\gamma_1 -b)/\gamma_1, & v<b\\
( v-\gamma_2 w + b^*\gamma_2 -b)/\gamma_2, & v\geq b
\end{cases} ,
\nonumber
\end{equation}
with $a/2 \leq b^* \leq (1-a)/2$ and $a/2 \leq b\leq (1+a)/2$.

%%%%%%%%%%%%%%%%%%%%%%%%%%%%%%%%%%%%%%%%%%%%%%%%%%%%%%%%%%%%%%%%%%%%%%%%%%%%%%%%%%%%%%%%%%%%%%%%%%%%%%%%%%%%%%%%%%
%%%%%%%%%%%%%%%%%%%%%%%%%%%%%%%%%%%%%%%%%%%%%%%%%%%%%%%%%%%%%%%%%%%%%%%%%%%%%%%%%%%%%%%%%%%%%%%%%%%%%%%%%%%%%%%%%%

The periodic orbit (as shown in Fig.~\ref{Fig:PWLML}) consists of four pieces labelled by $\mu=1,\ldots, 4$. On each section $\vec{z}_\mu(t) = \mathcal{E}_\mu(t) \vec{z}_\mu(0) + \mathcal{K}_\mu(t) \vec{c}_\mu$
where $\mathcal{E}_\mu(t) = \e^{A_\mu t}$, and $\mathcal{K}_\mu(t) = A_\mu^{-1}[\mathcal{E}_\mu(t)-I_m]$ with
$A_3=A_1$, $\vec{c}_3=\vec{c}_1$ and
\begin{align*}
A_1 & = \begin{bmatrix}  1/C & -1/C \\ 1/\gamma_2 & -1 \end{bmatrix}, & A_2 &= \begin{bmatrix}  -1/C & -1/C \\ 1/\gamma_2 & -1 \end{bmatrix},\\
A_4 & = \begin{bmatrix}  1/C & -1/C \\ 1/\gamma_1 & -1 \end{bmatrix} , & \vec{c}_1 &= \begin{bmatrix}  (I-a)/C \\ b^* - b/\gamma_2 \end{bmatrix}, \\
\vec{c}_2 &= \begin{bmatrix}  (1+I)/C \\ b^* - b/\gamma_2 \end{bmatrix}, & \vec{c}_4 &= \begin{bmatrix}  (I-a)/C \\ b^* - b/\gamma_1\end{bmatrix}.
\end{align*}
The periodic orbit crosses two switching planes with indicator functions $h_1(\vec{z}; 0) = v-b$ and $h_2(\vec{z}; 0) = v-(1+a)/2$. Therefore for initial data $\vec{z}(0) = (b, w(0))$, to find the periodic orbit we must solve the system of five nonlinear algebraic equations
\begin{align*} &v(T_1)= (1+a)/2, \qquad  v(T_1+T_2) = (1+a)/2 ,\\
&v(T_1+T_2 + T_3)=b,  \qquad v(\Delta) = b,  \qquad w(0) = w(\Delta) , \end{align*} for the switching times $T_1$, $T_2$, $T_3$, the period $\Delta$ and the initial value $w(0)$. 
This periodic orbit can be shown to be stable with Floquet exponent \cite{Coombes2008}
\[ r = -1 + (\Delta-2T_2)/(C\Delta).\]
Note that other types of periodic solutions are possible, such as those that cross only $v=b$ or $v=(1+a)/2$, or ones that do not cross either of the switching manifolds (as for a periodic orbit surrounding a linear center).

\begin{figure*}[t]
\begin{center}
\includegraphics[width=16cm]{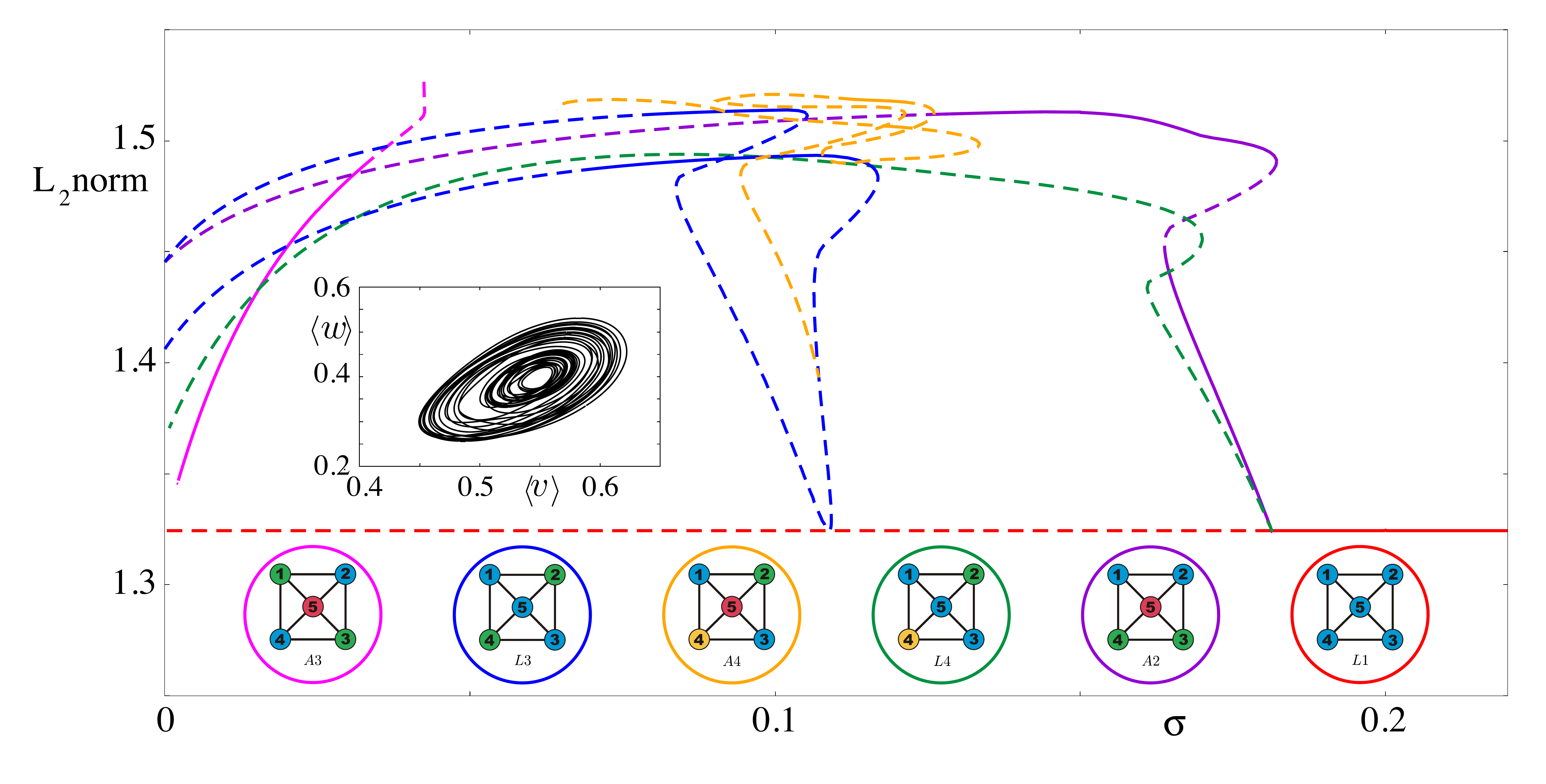}
\caption{Branches of cluster states and their stability upon variation of coupling strength $\sigma$ in the five-node network of nodes with PML local dynamics (parameter values as in Fig.~\ref{Fig:PWLML}) and diffusive coupling. Stable periodic solutions are denoted by solid lines while unstable solutions are denoted by dotted lines. Each branch is coloured according to the circle containing the relevant cluster state. In cases where conjugate solutions are also expected with identical stability properties (i.e. for the $A2$, $A4$, $L3$ and $L4$ branches) only one branch is shown. The inset shows the mean field dynamics for $\sigma= 0.065$. This kind of aperiodic behaviour dominates for $0.033 \lesssim \sigma \lesssim 0.071$. There is also a stable multirhythm for $0.0025 \lesssim \sigma \lesssim 0.0531$ (see main text).   }
\label{Fig:PMLBifDiag}
\end{center}
\end{figure*}

When nodes with PML local dynamics are diffusively coupled with the structure given by graph Laplacian \eqref{eq:Laplacian5} a number of different stable cluster states are observed. Similarly to the case of node dynamics described by the absolute model as discussed in section \ref{sec:Nonsmooth}, we can locate periodic orbits for cluster states and carry out linear stability calculations. We can then produce bifurcation diagrams such as that given in Fig.~\ref{Fig:PMLBifDiag}. Note that just as for the absolute model local dynamics, the PML model has a continuous vector field and therefore again all saltation matrices \eqref{eq:saltation} are the identity.
In practice, dealing with changes in the ordering of event crossings (defining a cluster) and to gain insight into what types of cluster state are possible (to provide good initial guesses for a numerical root-finding scheme) we also did simulations of a corresponding set of smooth ODEs, with a sufficiently steep switch to mimic the PML system.

%%%%%%%%%%%%%%%%%%%%%%%%%%%%%%%%%%%%%%%%%%%%%%%%%%%%%%%%%%%%%%%%%%%%%%%%%%%%%%%%%%%%%%%%%%%%%%%%%%%%%%%%%%%%%%%%%%
%%%%%%%%%%%%%%%%%%%%%%%%%%%%%%%%%%%%%%%%%%%%%%%%%%%%%%%%%%%%%%%%%%%%%%%%%%%%%%%%%%%%%%%%%%%%%%%%%%%%%%%%%%%%%%%%%%

The bifurcation diagram in Fig.~\ref{Fig:PWLML} shows the presence of $6$ of the possible $8$ cluster states. All but the $L4$ branch possess a range of values of $\sigma$ for which the periodic cluster state is stable. In the region $0.033 \lesssim \sigma \lesssim 0.071$ where there are no stable branches of periodic cluster states, aperiodic behaviour such as that shown in the inset of Fig.~\ref{Fig:PMLBifDiag} dominates. There is also a stable multirhythm,
where nodes can oscillate at different frequencies (not shown in Fig.~\ref{Fig:PMLBifDiag}). This has $\vec{x}_1(t) = \vec{x}_2(t+\Delta/2)$, $\vec{x}_3(t) = \vec{x}_4(t+\Delta/2)$ and $\vec{x}_5(t) = \vec{x}_5(t+\Delta/2)$ so $\vec{x}_5$ has a period of oscillation which is half that of $\vec{x}_i$, $i=1,\ldots,4$. The amplitudes of $\vec{x}_1$ and $\vec{x}_5$ are also very small resulting in a much smaller $L_2$-norm than the other solutions shown in Fig.~\ref{Fig:PMLBifDiag}.

\subsection{\label{sec:pwlif} Piecewise linear Integrate-and-fire network}

\begin{figure*}[t]
\begin{center}
\includegraphics[width=16cm]{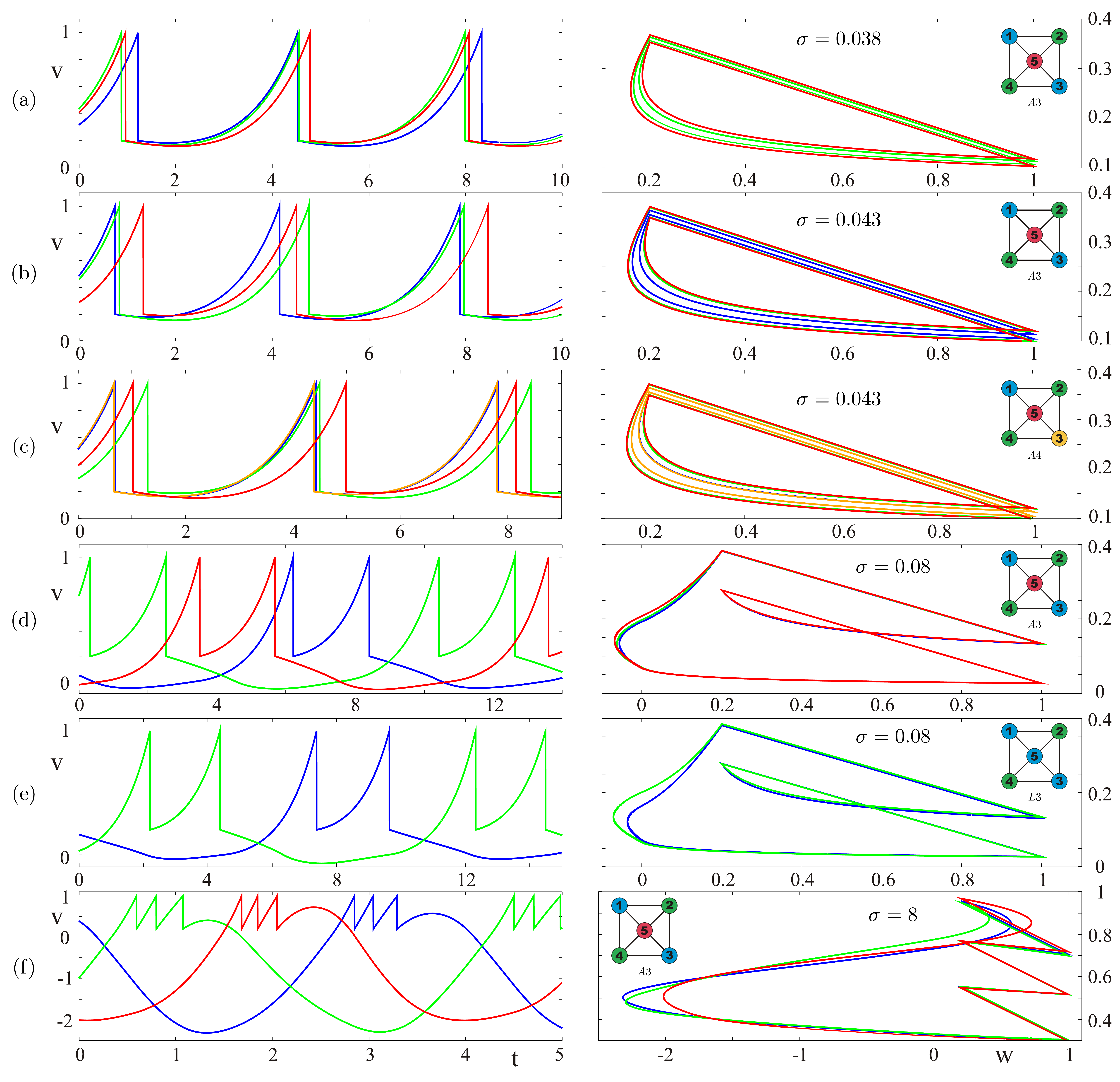}
\caption{Examples of cluster states for $\alpha=0.4$ and different values of coupling strength $\sigma$ in a network of five PWL-IF neurons. All parameters for the local dynamics are as in Fig.~\ref{Fig:pwlIF}, top left. The left panels show the time series for the voltage variable $v_i$ for each of the clusters within the state while the right panel shows each of the periodic orbits in the $(v_i,w_i)$-plane.
(a) A stable (A3) cluster state where each cluster fires twice per period observed for $\sigma=0.038$. In the right hand panel the orbit for $\vec{x}_1$ is obscured by that of $\vec{x}_5$ as (up to phase shift) they have virtually identical shape. (b) and (c) show the bistable $A3$ and $A4$ cluster states respectively for $\sigma=0.043$. As $\vec{x}_1$ and $\vec{x}_3$ are very close in the $A4$ cluster state, the orbit of $\vec{x}_3$ obscures that of $\vec{x}_1$.  (d) and (e) show the $A3$ and $L3$ cluster states which are bistable at $\sigma=0.08$. For large values of $\sigma$, the $A3$ cluster state remains stable but undergoes spike adding as in (e) where three spikes are observed within each cluster per orbit.  }
\label{Fig:PWLIFnetwork}
\end{center}
\end{figure*}

As our final example we consider a network of synaptically coupled PWL-IF neurons as discussed for the case of synchrony in section \ref{sec:MSF-PWL}. We continue to use the same five-node structured network \eqref{eq:Laplacian5}. We also choose for the local dynamics parameter values as in the top left of Fig.~\ref{Fig:pwlIF} that give tonic firing of each of the nodes when the coupling strength $\sigma$ is $0$. Setting the value of the synaptic rate parameter $\alpha=0.4$ we observe a number of different stable cluster states depending on the value of the coupling strength $\sigma$. A selection of these cluster states is given in
Fig.~\ref{Fig:PWLIFnetwork}. For these parameter values, the (tonic firing) synchronous state ($L1$) is stable for $0 \lesssim \sigma \lesssim 0.0334$. This undergoes a period-doubling bifurcation at $\sigma_{\rm{pd}} \simeq 0.0334$ to an $A3$ cluster state of doublets such as that pictured in Fig.~\ref{Fig:PWLIFnetwork} (a) for $\sigma \simeq 0.038$. For $0.0395 \lesssim \sigma \lesssim 0.0485$ this solution branch (as pictured for $\sigma \simeq 0.043$ in Fig.~\ref{Fig:PWLIFnetwork} (b)) is bistable with the periodic 4 cluster state $A4$ where $\vec{x}_1$ and $\vec{x}_3$ are not equal, but so close that it is hard to distinguish their time series and orbits (shown in blue and orange respectively in Fig.~\ref{Fig:PWLIFnetwork} (c) for $\sigma \simeq 0.043$). Further increasing $\sigma$ beyond 0.0485 leads to a loss of periodicity of solutions (both aperiodic $A4$ and $A5$ solution branches are observed). At around $\sigma=0.07$ again we see a stable periodic $A3$ cluster state for which all orbits now pass through the switch at $v_i=0$. This is pictured in Fig.~\ref{Fig:PWLIFnetwork} (d) for $\sigma \simeq  0.08$ where the solution is bistable with the two cluster $L3$ solution as in Fig.~\ref{Fig:PWLIFnetwork} (e). The $A3$ cluster state remains stable for $\sigma \gtrsim 0.07$ undergoing a series of spike adding bifurcations for large values of $\sigma$. Fig.~\ref{Fig:PWLIFnetwork} (f) shows the solution branch when $\sigma=8$ where each cluster spikes three times per orbit.  We note that recent work by Chen \textit{et al}. \cite{Chen2017} has also considered cluster states in leaky IF networks with $\alpha$-function synapses, although their results are restricted to weak-coupling in small all-to-all coupled networks.

%%%%%%%%%%%%%%%%%%%%%%%%%%%%%%%%%%%%%%%%%%%%%%%%%%%%%%%%%%%%%%%%%%%%%%%%%%%%%%%%%%%%%%%%%%%%%%%%%%%%%%%%%%%%%%%%%%
%%%%%%%%%%%%%%%%%%%%%%%%%%%%%%%%%%%%%%%%%%%%%%%%%%%%%%%%%%%%%%%%%%%%%%%%%%%%%%%%%%%%%%%%%%%%%%%%%%%%%%%%%%%%%%%%%%

\section{\label{sec:Discussion}Discussion}

In this paper we have shown how to extend the well known MSF formalism for networks of smooth dynamical systems to treat PWL models relevant to the modelling of coupled neural oscillators.  Moreover, we have also augmented the recent extension of the MSF approach that treats synchronous cluster states (exploiting spatial or pointwise network symmetries and balanced colourings).  We have illustrated the general theory using some relatively simple networks, though with various choices of node dynamics.  This has highlighted that the degree of synchronisation or type of clustering pattern observed in neuronal networks may have as much to do with the choice of connectivity as it has  with the choice of node dynamics.  This reinforces the notion that care must be taken when trying to recover structural connectivity from functional connectivity \cite{Stroud2015,Alderisio2017}.  Similarly it is known that the reconstruction of network dynamics from observables can be strongly influenced by symmetries \cite{Whalen2015}.  Although we have focused on example systems drawn from theoretical neuroscience, there are many PWL models in physics and engineering for which the framework we have presented can also be applied, such as electromechanical oscillator networks \cite{Shiroky2016}.
Interestingly, recent work by Cho \textit{et al}. \cite{Cho2017} has made a connection between synchronised cluster states and \textit{chimeras} \cite{kuramoto-battogtokh-2002,Abrams2004}, where a sub-population of oscillators synchronises in an otherwise incoherent sea.  By exploiting the symmetries of the network they show how to construct cluster partitions, in which the synchronisation stability of the individual clusters in each set is decoupled from that in all the other sets.  Using this they show that a chimera can exist if the fully synchronous solution is unstable and there is a cluster state with a main stable group, and all other clusters unstable.  Thus it would be of interest to revisit this work, formulated in a smooth setting, using the PWL perspective presented here.

It is well to mention that we have only considered networks of \textit{identical} oscillators, and it is well known that heterogeneity can disrupt synchrony \cite{Denker2002}.  The MSF extension to nearly identical units is relatively straightforward \cite{Sun2009}, though for PWL systems there is also the further possibility of the explicit construction of network states (by patching together local trajectories).  However, the MSF is not the only tool for analysing cluster states and it would be interesting to consider the augmentation of Lyapunov function methods, such as those developed by Belyk \textit{et al}. \cite{Belykh2003}, to treat nonsnmooth systems with event driven interactions.  In a neural context it would also be important to treat delays (arising from the finite speed of action potential propagation), and there is a growing body of work in this area that could be revisited from a PWL perspective \cite{Flunkert2010,Dahms2012,Huddy2016}.  Indeed delays are well known to lead to oscillatory instabilities, which at the network level may manifest as travelling waves.  This highlights an important mathematical extension of the work presented here, namely that periodic states have phase shift symmetries which must also be accounted for. When combined with the pointwise network symmetries these lead to spatio-temporal symmetries \cite{Buono2001, Golubitsky2016}.

A periodic solution $\vec{x}(t)$ has a group $K$ of spatial symmetries (the $\gamma \in \Gamma$ such that $\gamma \vec{x}(t) = \vec{x}(t)$ for all $t$). It is these symmetries which we have discussed in this paper. However, periodic orbits also possess a group $H$ of spatial symmetries under which the periodic orbit is invariant (that is the $\gamma \in \Gamma$ such that $\{\gamma \vec{x}(t) \ : \ t \in \mathbb{R}\} =\{\vec{x}(t) \ : \ t \in \mathbb{R}\}$). Thus there exists a phase shift $\phi \in \mathbb{S}^1$ such that $\gamma \vec{x}(0) = \vec{x}(\phi)$ and therefore $\gamma \vec{x}(t-\phi) = \vec{x}(t)$ for all $t$. The pair $(\gamma, \phi)$ define a spatio-temporal symmetry of the periodic orbit. The group $K$ is normal in $H$ and $H/K$ is a finite subgroup of $\mathbb{S}^1$. If $H=K$ then there are no phase-shift symmetries. The pair $(H,K)$ determines the pattern of phase shifts (up to isomorphism). Patterns  of  phase  shifts  which may exist for a given network structure with symmetry group $\Gamma$ can be determined using the $H/K$ theorem \cite{Buono2001}. This theorem gives conditions on subgroups $K \subseteq H \subseteq \Gamma$ such that the pair $(H,K)$ could determine the spatio-temporal symmetries of a hyperbolic periodic solution of the $\Gamma$-equivariant network equations, allowing one to compute a catalogue of possible phase shifted states. These may include so called multirhythm states (such as that found here for the PWL Morris--Lecar local dynamics on the five-node network, see section \ref{sec:pwlml}) where some cells are forced to oscillate with a frequency that is a rational multiple of that of the other cells \cite{Golubitsky2016}. Other phase shifted cluster states may also arise due to additional constraints imposed by the choice of Laplacian or balanced coupling. There are also approaches exploiting spatio-temporal symmetries for computing stability of periodic states near Hopf bifurcation \cite{Golubitsky1988} which may also be extendable to phase-shifted cluster states from Laplacian or balanced coupling.  The extension of our work to treat spatio-temporal symmetries will be presented elsewhere.

\begin{acknowledgments}
The authors would like to thank Yi Ming Lai, Mason Porter, and R\"{u}diger Thul for helpful comments made on a first draft of this paper.
This work was supported by the Engineering and Physical Sciences Research Council [grant number EP/P007031/1].
\end{acknowledgments}

\appendix
\section{\label{sec:saltation}Saltation matrices}

%%%%%%%%%%%%%%%%%%%%%%%%%%%%%%%%%%%%%%%%%%%%%%%%%%%%%%%%%%%%%%%%%%%%%%%%%%%%%%%%%%%%%%%%%%%%%%%%%%%%%%%%%%%%%%%%%%
\begin{figure}[ht]
\begin{center}
\includegraphics[width=8cm]{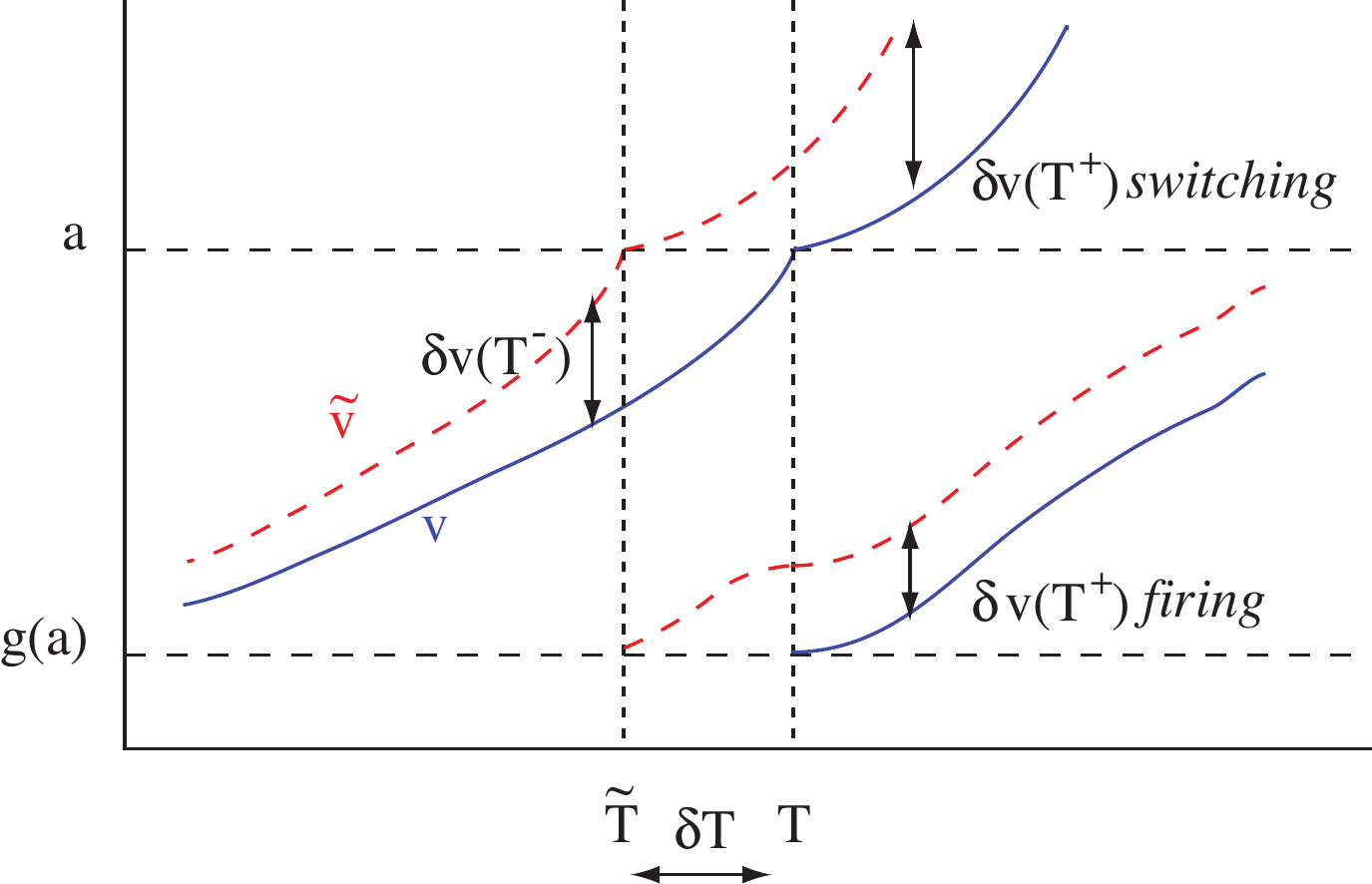}
\caption{
Evolution of a perturbation at a \textit{switching} or \textit{firing event}. The solid blue line is the trajectory of an unperturbed orbit, with the event at time $T$.  The dashed red line is the perturbed trajectory with an event at time $\widetilde{T}$.  In this illustration the deviation in event times is ${\delta}T = \widetilde{T}-T<0$.
}
\label{Fig:dtneg}
\end{center}
\end{figure}
%%%%%%%%%%%%%%%%%%%%%%%%%%%%%%%%%%%%%%%%%%%%%%%%%%%%%%%%%%%%%%%%%%%%%%%%%%%%%%%%%%%%%%%%%%%%%%%%%%%%%%%%%%%%%%%%%%
Saltation matrices are a useful way to augment the standard theory of smooth dynamical systems to treat the evolution of perturbations through regions in phase phase where there is a discontinuity in the vector field or the flow.  They have recently been used in a PWL network setting in \cite{Coombes2016}, and for nonlinear IF networks in \cite{Ladenbauer2013}.  Here we combine these two perspectives to construct the saltation matrices for hybrid systems that may have switches in their dynamics as well as discontinuities arising from reset.  Given that these cannot occur at the same time we adopt a notation to track either of these types using an indicator function $h(\vec{z})$ such that an event (switching or reset) occurs when $h(\vec{z}(T))=0$ (with different choices of $h$ for switching and reset).
The state of the system after an event is given by $\vec{z}(T^+)=\vec{g}(\vec{z}(T))$. Let us now consider an arbitrary trajectory of a dynamical system $\dot{\vec{z}} = \vec{F}(\vec{z})$ and linearise around this trajectory.  Between switching or firing events small perturbations $\delta \vec{z}$ satisfy
\begin{equation}
\FD{}{t} \delta \vec{z}(t) = D \vec{F}(\vec{z}(t)) \delta \vec{z}(t), \quad \delta \vec{z}(0) = \delta \vec{z}_0 .
\nonumber
\end{equation}
For the PWL models discussed in this paper then $D \vec{F}(\vec{z}(t))$ is a piecewise constant matrix, so that
$\delta \vec{z}(t) = \exp(A t) \delta \vec{z}_0$ with $A$ a constant matrix (such as $A_{L,R}$ in \S \ref{sec:MSF-PWL}).
A perturbed trajectory defined by $\widetilde{\vec{z}}(t)$ can be constructed as $\widetilde{\vec{z}}(t) = \vec{z}(t)+\delta \vec{z}(t)$. We can then define a corresponding perturbed event time $\widetilde{T}$ according to the condition $h(\widetilde{\vec{z}}(\widetilde{T}))=0$.  We shall denote the difference in the unperturbed and perturbed event times as $\delta T = \widetilde{T}-T$. For the case that $\delta T <0$ we have that
$\widetilde{\vec{z}}(T^+)=\widetilde{\vec{z}}(\widetilde{T}^+-\delta T) \simeq \vec{z}(\widetilde{T}^+) - \dot{\vec{z}}(\widetilde{T}^+) \delta T$, and see Fig.~\ref{Fig:dtneg}.
Hence,
\begin{align}
&\widetilde{\vec{z}}({T}^{+}) \simeq \vec{g}(\vec{z}(\widetilde{T}^{-})) - \dot{\vec{z}}(\widetilde{T}^{+}){\delta}T \nonumber \\
&\simeq \vec{g}(\widetilde{\vec{z}}({T}^{-})+\dot{\widetilde{\vec{z}}}({T}^{-}){\delta}T) - \dot{\vec{z}}(\widetilde{T}^{+}){\delta}T \nonumber \\
&\simeq \vec{g}({\vec{z}}({T}^{-})+\delta \vec{z}({T}^{-})+\dot{\widetilde{\vec{z}}}({T}^{-}){\delta}T) - \dot{\vec{z}}(\widetilde{T}^{+}){\delta}T .
\nonumber
\end{align}
By performing a Taylor expansion of $\vec{g}$ we obtain:
\begin{align}
\widetilde{\vec{z}}({T}^{+}) &\simeq \vec{z}({T}^{+}) + D \vec{g} (\vec{z}({T}^{-}))
\left[\delta \vec{z}({T}^{-})+\dot{\widetilde{\vec{z}}}({T}^{-}){\delta}T\right] \nonumber \\
&- \dot{\vec{z}}(\widetilde{T}^{+}){\delta}T ,
\label{spam}
\end{align}
where we have used the result that $\vec{g}(\vec{z}({T}^{-}))=\vec{z}({T}^{+})$.  To calculate $\delta T$ we Taylor expand the indicator function for the perturbed trajectory as
\begin{align}
&h(\widetilde{\vec{z}}(\widetilde{T}^{-})) = h(\vec{z}({T}^{-}+\delta T) + \delta \vec{z} ({T}^{-}+\delta T) ) \nonumber \\
&\simeq h(\vec{z}({T}^{-})+ \dot{\vec{z}} (T^-)\delta T + \delta \vec{z} ({T}^-) ) \nonumber \\
&\simeq h(\vec{z}({T}^{-})) + \nabla_{\vec{z}} h(\vec{z}({T}^{-})) \cdot [ \dot{\vec{z}} (T^-)\delta T + \delta \vec{z} ({T}^-)].
\nonumber
\end{align}
Using the fact that $h(\widetilde{\vec{z}}(\widetilde{T}^{-})) = 0=h({\vec{z}}({T}^{-}))$ we have that
\begin{equation}
\delta T = -  \frac{\nabla_{\vec{z}} h(\vec{z}({T}^{-})) \cdot \delta \vec{z} ({T}^-)}{\nabla_{\vec{z}} h(\vec{z}({T}^{-})) \cdot \dot{\vec{z}} (T^-)} .
\label{deltaT}
\end{equation}
Now making use of the approximations $\dot{\widetilde{\vec{z}}}(T^-) \simeq \dot{{\vec{z}}}(T^-)$ and $\dot{{\vec{z}}}(\widetilde{T}^+) \simeq \dot{{\vec{z}}}({T}^+)$ (and see Fig.~\ref{Fig:dtneg}), we may use (\ref{spam}) and (\ref{deltaT}) to construct $\delta \vec{z}(T^+)$ in the form
\begin{equation}
\delta \vec{z}(T^+) \simeq K(T) \delta \vec{z}(T^-),
\nonumber
\end{equation}
where $K(T)$ is the saltation matrix:
\begin{align}
&K(T) = D \vec{g}(\vec{z}(T^-)) \nonumber \\
& +   \frac{\left [\dot{\vec{z}}(T^+) - D \vec{g}(\vec{z}(T^-)) \dot{\vec{z}}(T^-) \right ] \left [ \nabla_{\vec{z}} h(\vec{z}(T^-)) \right ]^\mathsf{T}}{\nabla_{\vec{z}} h(\vec{z}(T^-)) \cdot \dot{\vec{z}} (T^-)} .
\nonumber
\end{align}
The saltation matrix allows us to relate perturbations just before an event to those just after.  A similar calculation with $\delta T >0$ gives precisely the same formula for $K(T)$.  We note that for a switching event (where reset does not occur) then $\vec{g}(\vec{z})=\vec{z}$ and $D \vec{g}(\vec{z})$ is the identity matrix.
Moreover, if in this case the vector field is continuous then $\dot{\vec{z}}(T^+)=\dot{\vec{z}}(T^-)$ and $K(T)$ reduces to the identity matrix as expected.

For the single PWL-IF node $\vec{z} \in \Rset^2$ with $\nabla_{\vec{z}} h(\vec{z};v_\text{th})=(1,0)^\mathsf{T}$, and
\begin{equation}
D \vec{g}(\vec{z}) = \begin{bmatrix}
0 & 0 \\
0 & 1
\end{bmatrix}.
\nonumber
\end{equation}
For the PWL-IF network  $\vec{z}_i \in \Rset^4$ with $\nabla_{\vec{z}} h(\vec{z}_i;v_\text{th})=(1,0,0,0)^\mathsf{T}$, and
\begin{equation}
D \vec{g}(\vec{z}_i) = \begin{bmatrix}
0 & 0 & 0 & 0\\
0 & 1 & 0 & 0\\
0 & 0 & 1 & 0\\
0 & 0 & 0 & 1
\end{bmatrix}.
\nonumber
\end{equation}

%\bibliography{Lucie}

%

\end{document}